\documentclass[journal]{IEEEtran}

\newtheorem{thm}{Theorem}
\newtheorem{defn}[thm]{Definition}
\newtheorem{lem}[thm]{Lemma}
\newtheorem{col}[thm]{Corollary}
\newtheorem{ex}[thm]{Example}

\newtheorem{remark}[thm]{Remark}
\newtheorem{prp}[thm]{Proposition}
\ifCLASSINFOpdf
\else
\fi
\usepackage{tikz}
\usepackage{tkz-graph}
\usepackage{graphicx}      

\usepackage{algorithm,algorithmic}
\usepackage{amsmath,amssymb,graphicx}
\usepackage{amsfonts}
\usepackage{mathtools}
\usepackage{caption}
\usepackage{subcaption}

\DeclareMathOperator{\sign}{sign}
\DeclareMathOperator{\sgn}{sgn}
\DeclareMathOperator{\spec}{spec}

\usepackage{tcolorbox}

\usepackage[colorlinks=true,linkcolor=green]{hyperref}
\hypersetup{
	colorlinks,
	linkcolor={blue!80!black},
	citecolor={green!80!black},
	urlcolor={brown!80!black}
}
\def\equationautorefname~#1\null{Equation~(#1)\null}
\usepackage[prependcaption,colorinlistoftodos]{todonotes}

\begin{document}
%
\title{Strong Structural Controllability of Systems on Colored Graphs}


\author{\IEEEauthorblockN{Jiajia Jia,
		Harry L. Trentelman,
		Wouter Baar,
		and
		Kanat M. Camlibel}

\thanks{The authors are with the Bernoulli Institute for Mathematics, Computer Science and Artificial Intelligence, University of Groningen, The Netherlands (e-mail: j.jia@rug.nl;h.l.trentelman@rug.nl;w.baar.1@student.rug.nl;
	m.k.camlibel@rug.nl)}

\markboth{Draft submission: Transactions of Automatic Control,~Vol.~XX, No.~XX, Month~year}%
{Shell \MakeLowercase{\textit{et al.}}:
	Draft submission: Transactions of Automatic Control,~Vol.~XX, No.~XX, Month~year}}


%



\IEEEtitleabstractindextext{%
%
%
%
%
%

\begin{abstract}
This paper deals with structural controllability of leader-follower networks. The system matrix defining the network dynamics is a pattern matrix in which a priori given entries are equal to zero, while the remaining entries take nonzero values. The network is called strongly structurally controllable if for all choices of real values for the nonzero entries in the pattern matrix, the system is controllable in the classical sense. In this paper we introduce a more general notion of strong structural controllability which deals with the situation that given nonzero entries in the system's pattern matrix are constrained to take identical nonzero values.
The constraint of identical nonzero entries can be caused by symmetry considerations or physical constraints on the network. The aim of this paper is to establish graph theoretic conditions for this more general property of strong structural controllability.
\end{abstract}

}

\maketitle

\IEEEdisplaynontitleabstractindextext

%
\IEEEpeerreviewmaketitle

\section{Introduction}

The past two decades have shown an increasing research effort in networked dynamical systems. To a large extend this increase has been caused by technological developments such as the emergence of the internet and the growing relevance of smart power grids. Also the spreading interest in social networks and biological systems have contributed to this surge \cite{Newman2003,ME2010,liu2016control,Ruths2014}. 

A fundamental issue in networked systems is that of controllability. This issue deals with the question whether
all parts of the global network can be adequately influenced or manipulated by applying control inputs only locally to the network. A vast amount of literature has been devoted to several variations on this issue, see \cite{tanner2004,LCWX2008,RJME2009,LYB2011,EMCCB2012,WCT2017} and the references therein. In most of the literature, a networked  system is a collection of input-state-output systems, called agents, together with an interconnection structure between them. Some of these systems can also receive input from outside the network, and are called leaders. The remaining systems are called followers. At a higher level of abstraction, a networked system can be described by a direct graph, called the network graph, where the vertices represent the input-state-output systems and the edges represent the interactions between them. Controllability of the networked system then deals with the question whether the states of all agents can be steered from any initial state to any final state in finite time by applying suitable input signals  to the network through the leaders.

Using the observation that the underlying graph plays an essential role in the controllability properties of the networked system \cite{EMCCB2012}, an increasing amount of literature has been devoted to uncovering this connection, see  \cite{GP2013,MM2013,CM2014} and the references therein. In order to be able to zoom in on the role of the network graph, it is common to proceed with the simplest possible dynamics at the vertices of the graph, and to take the agents to be single integrators, with a one-dimensional state space. These single integrators are interconnected through the network graph, and the interconnection strengths are given by the weights on the edges. Based on this, the overal networked system can obviously be represented by a linear input-state-output system of the form 
\[
\dot{x} = Ax + Bu,
\]
where the system matrix $A \in \mathbb{R}^{n \times n}$ represents the network structure with the given edge weights, and the matrix $B \in \mathbb{R}^{n \times m}$ encodes which $m$ vertices are the leaders. The $n$-dimensional state vector $x$ consists of the states of the $n$ agents, and the $m$-dimensional vecor $u$ collects the input signals to the $m$ leader vertices.

Roughly speaking, the research on network controllability based on the above model can be subdivided into three directions. The first direction is based on the assumption that the edge weights in the network are known exactly. In this case the matrix $A$ is a given constant matrix, and specific dynamics is considered for the network. For example, the system matrix can be defined as the adjacency matrix of the graph   \cite{CS2010}, or the graph Laplacian matrix \cite{tanner2004,RJME2009,EMCCB2012,MH2016,ZCC2014,ZCC2011}. Furthermore, a framework for controllability was also introduced in \cite{YZDW2013}, offering tools to treat controllability of complex networks with arbitrary structure and edge weights. Related results can be found in \cite{YZ2014,nie2015}. We also refer to \cite{TT2018,TR2014}.

A second research direction deals with the situation that the exact values of the edge weights are not known, but only information on whether these weights are zero or nonzero is available. In this case, the system matrix is not a known, given, matrix, but rather a matrix with a certain zero/nonzero pattern: some of the entries are known to be equal to zero, the other entries are unknown. This framework deals with the concept of {\em structural controllability}. Up to now, two types of structural controllability have been studied, namely {\em weak} structural controllability and {\em strong} structural controllability. A networked system of the form above is called weakly structurally controllable if there exists at least one choice of values for the nonzero entries in the system matrices such that the corresponding matrix pair $(A,B)$ is controllable. The networked system is called strongly structurally controllable if, roughly speaking, for {\em all} choices of values for the nonzero entries the matrix pair $(A,B)$ is controllable.
Conditions for weak and strong structural controllability can be expressed entirely in terms of the underlying network graph, using concepts like cactus graphs, maximal matchings, and zero forcing sets, see \cite{Lin74,MY1979,LYB2011,CM2013,MZC2014,WRS2014,TD2015}.

A third, more recent, research direction again deals with weak and strong structural controllability. However, the nonzero entries in the pattern matrices defining the networked system can no longer take arbitrary nonzero real values, independently of each other. Instead, in this framework the situation is considered that there are certain constraints on some of the nonzero entries. These constraints can require that some of the nonzero entries have given values, see e.g.  \cite{MHM2018}, or that there are given linear dependencies between some of the nonzero entries, see \cite{LM2017}. In particular, in \cite{LM2017} necessary and sufficient conditions for weak structural controllability were established in terms of {\em colored graphs.}

The present paper contributes to this third research direction. In the context of strong structural controllability it deals with the situation that the nonzero entries in the system matrix can no longer take arbitrary values. Instead, the values of certain a priori specified nonzero entries in the system matrix are constrained to be identical. 
Obviously, in real world networks it is indeed a typical situation that certain edge weights are equal, either by symmetry considerations or by the physics of the underlying problem.
An example is provided by the case of undirected networks, in which the network graph has to be symmetric. Another example is provided by networked systems defined in terms of  so-called network-of-networks \cite{CAM2014}, which are obtained by taking the Cartesian product of smaller factor networks. 
For each factor network, the internal edge weights are independent. However, by applying the Cartesian product, some edge weights in the overal network will become identical.
In the present paper we will establish conditions for this new notion of strong structural controllability in terms of colored graphs.

The main contributions of this paper are as follows:
\begin{enumerate}
	\item We introduce a new color change rule and define the corresponding notion of zero forcing set. To do this, we consider colored bipartite graphs and establish a necessary and sufficient graph-theoretic condition for nonsingularity of the pattern class associated with this bipartite graph.
	\item We provide a sufficient graph theoretic condition for our new notion of strong structural controllability in terms of zero forcing sets.
	\item We introduce so called elementary edge operations that can be applied to the original network graph and that preserve the property of strong structural controllability.
	\item A sufficient graph theoretic condition for strong structural controllability is developed based on the notion of edge-operations-color-change derived set which is obtained by applying elementary edge operations and the color change rule iteratively.
\end{enumerate}

The organization of this paper is as follows. In Section 2,  some preliminaries are presented. 
In Section 3,  we give a formal definition of the main problem treated in this paper in terms of systems defined on colored graphs. 
In Section 4, we establish our main result, giving a sufficient graph-theoretic condition for strong structural controllability of systems defined on colored graphs.
Section 5 provides two additional sufficient graph-theoretic conditions. To establish these conditions, we introduce the concept of elementary edge operations and the associated notion of edge-operations-color-change derived set. This set is obtained from the initial coloring set by an iterative procedure involving successive and alternating applications of elementary edge operations and the color change rule.   
Finally,  Section 6 formulates the conclusions of this paper.
We note that a preliminary version \cite{JTBC2018} of this paper has appeared in the proceedings of NecSys 2018. In that note, the condition for strong structural controllability in terms of our new concept of zero forcing set was stated without giving any of the proofs. The present paper provides these proofs, and in addition provides new conditions for strong structural controllability in terms of elementary edge operations and the concept of edge-operations-color-change derived set that were not yet given in \cite{JTBC2018}.

\section{Preliminaries} 

In this paper, we will use standard notation. Let $\mathbb{C}$ and $\mathbb{R}$ denote the fields of complex and real numbers, respectively. 
The spaces of $n$-dimensional real and complex vectors are denoted by $\mathbb{R}^{n}$ and 
$\mathbb{C}^{n}$, respectively. Likewise, the spaces of $n \times m$ real and complex matrices are denoted by $\mathbb{R}^{n \times m}$ and $\mathbb{C}^{n \times m}$, respectively.
For a given $n \times m$ matrix $A$, the entry in the $i$th row and $j$th column is denoted by $A_{ij}$.  
For a given $m \times n$ matrix $A$ and for given subsets 
$S = \{s_1,s_2,\ldots,s_k\} \subseteq \{1,2,\ldots,m\}$ and $T = \{t_1,t_2,\ldots,t_l\} \subseteq \{1,2,\ldots,n\}$ we define the $k \times l$ submatrix of $A$ associated with $S$ and $T$ as the matrix $A_{S,T}$ with $(A_{S,T})_{ij} := A_{s_i t_j}$. Similarly, for a given $n$-dimensional vector $x$, we denote by $x_{T}$ the subvector of $x$ consisting of the entries of $x$ corresponding to $T$.
For a given square matrix $A$, we denote its determinant by $\det(A)$. Finally, $I$ and $\mathbf{0}$ will denote the identity and zero matrix of appropriate dimensions, respectively.

\subsection{Elements of Graph Theory}

Let $\mathcal{G} = (V,E)$ be a directed graph, with vertex set $V = \{1,2, ...,n\}$,  and the edge set $E$ a subset of  $V \times V$.
In this paper, we will only consider simple graphs, that is, the edge set $E$ does not contain edges of the form $(i,i)$. In our paper, the phrase `directed graph' will always refer to a simple directed graph.
We call vertex $j$ an out-neighbor of vertex $i$ if $(i,j) \in E$. 
We denote by  $N(i) := \{j \in V \mid (i,j) \in E \}$ the set of all out-neighbors of $i$.
 Given a subset $S$ of the vertex set $V$ and a subset $X \subseteq S$, we denote by
\[
N_{V \setminus S}(X) = \{j \in V \setminus S \mid \exists~ i \in X ~\mbox{such that} ~(i,j) \in E \},
\]
the set of all vertices outside $S$, but an out-neighbor of some vertex in $X$.   
A  directed graph $\mathcal{G}_{1} = (V_{1},E_{1})$ is called a subgraph of $\mathcal{G}$ if $V_{1} \subseteq V$ and $E_{1} \subseteq E$.

Associated with a given directed graph $\mathcal{G} = (V,E)$ we consider the set of matrices
\[
\mathcal{W}(\mathcal{G}) := \{ W \in \mathbb{R}^{n \times n} \mid W_{ij} \neq 0 \mbox{ iff } (j,i) \in E \}.
\]
For any such $W$ and $(j,i) \in E$, the entry $W_{ij}$ is called the weight of the edge $(j,i)$.
Any such matrix $W$ is called a {\em weighted adjacency matrix} of the graph.
For a given directed graph $\mathcal{G} = (V,E)$, we denote the associated graph with weighted adjacency matrix $W$ by $\mathcal{G}(W) = (V,E,W)$. 
This is then called the {\em  weighted graph} associated with the graph $\mathcal{G}= (V,E)$ and weighted adjacency matrix $W$. 
Finally, we define the graph  $\mathcal{G} = (V,E)$ to be an {\em undirected graph} if $(i,j) \in E$ whenever $(j,i) \in E$. 
In that case the order of $i$ and $j$ in $(i,j)$ does not matter and we interpret the edge set $E$ as the set of unordered pairs $\{i,j\}$ where  $(i,j) \in E$.

An undirected graph $\mathcal{G} = (V,E)$ is called {\em bipartite} if there exist nonempty disjoint subsets $X$ and $Y$ of $V$ such that $X \cup Y = V$ and $\{i,j\} \in E$ only if $i \in X$ and $j \in Y$. 
Such {\em bipartite graph} is denoted by  $G = (X,Y,E_{XY})$ where we denote the edge set by $E_{XY}$ to stress that it contains edges $\{i,j\}$ with $i \in X$ and $j \in Y$. In this paper we will use the symbol $\mathcal{G}$ for arbitrary directed graphs and $G$ for bipartite graphs.

A set of $t$ edges $  {m} \subseteq E_{XY}$  is called a {\em $t$-matching} in $G$, if no two distinct edges in $  {m}$ share  a vertex.  
In the special case that $|X| = |Y| =t$, such a $t$-matching is called a {\em perfect matching}.

For a bipartite graph $G = (X,Y,E_{XY})$, with vertex sets $X$ and $Y$ given by $X = \{x_{1},x_{2}, \ldots, x_{s}\}$ and $Y = \{y_{1},y_{2}, \ldots, y_{t}\}$,
we define {\em the pattern class} of $G$ by
\[\mathcal{P}(G) =
\{M \in \mathbb{C}^{t \times s}\mid M_{ji} \neq 0 \mbox{ iff } \{x_{i},y_{j}\} \in E_{XY}\}.\]
Note that matrices $M \in \mathcal{P}(G)$ may not be square since the cardinalities of $X$ and $Y$ can differ. 
Also note that, in the context of pattern classes for undirected bipartite graphs, we allow complex matrices.

\subsection{Controllability of Systems Defined on  Graphs}

For a  directed graph $\mathcal{G} = (V,E)$ with vertex set $V = \{1,2,\ldots,n\}$, the {\em qualitative class} of $\mathcal{G}$ is defined as the family of matrices 
\[
\mathcal{Q}(\mathcal{G})  = \{A \in \mathbb{R}^{n \times n} \mid \mbox{for } i \neq j: A_{ij}
\neq 0 \mbox{ iff } (j,i) \in E \}.
\]
Note that the diagonal entries of $A \in \mathcal{Q}(\mathcal{G})$ do not depend on the structure of $\mathcal{G}$ and can take arbitrary real values. 

Next, we specify a subset $V_{L}  = \{v_{1},v_2, \ldots, v_{m}\}$ of $V$, called the  {\em the leader set}, and consider the following family of leader/follower systems defined on the graph $\mathcal{G}$ with dynamics
\begin{equation}\label{e:system}
\dot{x} = Ax + Bu ,
\end{equation}
where $x \in \mathbb{R}^{n}$ is the state and $u \in \mathbb{R}^{m}$ is the input.
The systems \eqref{e:system} have the distinguishing feature that the matrix $A$ belongs to $\mathcal{Q}(\mathcal{G})$ and $B = B(V;V_L)$ 
is defined as the $n \times m$ matrix given by
\begin{equation}\label{e:inputmatrix}
B_{ij} = \left\{ \begin{array}{lcl} 1 ~ \mbox{if} ~ i = v_{j},\\ 0 ~\mbox{otherwise}.\end{array}\right.
\end{equation}

An important notion   associated with systems defined on a graph $\mathcal{G}$  as in \eqref{e:system} is the notion of strong structural controllability.
\begin{defn}\label{d:cos}
	Let $\mathcal{Q}' \subseteq \mathcal{Q}(\mathcal{G})$.
	The system defined on the directed graph $\mathcal{G}=(V,E)$  with dynamics  \eqref{e:system}  and leader set $V_{L} \subseteq V$ is called   {\em strongly structurally controllable with respect to $\mathcal{Q}'$} if the pair $(A,B)$ is controllable for all $A \in \mathcal{Q}'$. 
	In that case we will simply say that  $(\mathcal{G};V_{L})$ is controllable with respect to $\mathcal{Q}'$. 
\end{defn}

One special case of the above notion is that $(\mathcal{G};V_{L})$ is controllable with respect to $\mathcal{Q}(\mathcal{G})$. 
In that case, we will simply say that $(\mathcal{G};V_{L})$ is controllable.
Another special case is that $(\mathcal{G};V_{L})$ is controllable with respect to  $\mathcal{Q}' $ where, for a given weighted adjacency matrix $W \in \mathcal{W}(\mathcal{G})$, $\mathcal{Q}'$ is the subclass of $\mathcal{Q}(\mathcal{G})$ defined by
\[\mathcal{Q}_{W}(\mathcal{G})  = \{A \in \mathcal{Q}(\mathcal{G}) \mid \mbox{for}~ i \neq j: A_{ij} = W_{ij} \}.\]
This subclass is called the {\em weighted qualitative class} associated with $W$.
Note that the off-diagonal elements of $A \in \mathcal{Q}_{W}(\mathcal{G})$ are fixed by those of the given adjacency matrix, while, again, the diagonal entries of $A \in \mathcal{Q}_{W}(\mathcal{G})$ can take arbitrary real values. Obviously
\[\mathcal{Q}(\mathcal{G}) = \bigcup_{W \in \mathcal{W}(\mathcal{G})} \mathcal{Q}_{W}(\mathcal{G}).\]
Since there is a unique weighted graph $\mathcal{G}(W) = (V,E,W)$ associated with the graph $\mathcal{G}= (V,E)$ and weighted adjacency matrix $W$, we will simply say that {\em $(\mathcal{G}(W);V_{L})$ is controllable} if  $(\mathcal{G};V_{L})$ is controllable with respect to $\mathcal{Q}_{W}(\mathcal{G})$. 

\subsection{Zero Forcing Set and Controllability of $(\mathcal{G};V_{L})$}\label{s:ZFS}

Let $\mathcal{G} = (V,E)$ be a  directed graph with vertices colored either black or white.  
We now introduce the following color change rule \cite{A2008}:
if $v$ is a black vertex in $\mathcal{G}$ with exactly one white
out-neighbor $u$, then we change the color of $u$ to black, and write $v \xrightarrow{c} u$.  
Such a color change is called a {\em force}.
A subset $C$ of $V$ is called a {\em coloring set} if the vertices in $C$ are initially colored black and those in $V \setminus C$ initially colored white.
Given a coloring set $C \subseteq V$,  the derived set $\mathcal{D}(C)$ is the set of black vertices obtained after repeated application of the color change rule, until no more changes are possible.
It was shown in \cite{A2008} that the derived set is indeed uniquely defined, in the sense that it does not depend on the order in which the color changes are applied to the original coloring set $C$.
A coloring set $C \subseteq V$ is  called a {\em zero forcing set for} $\mathcal{G}$ if $\mathcal{D}(C) = V$.
Given a zero forcing set for $\mathcal{G}$, we can list the  forces in the order in which they were performed to color  all vertices in the graph black. 
Such a list is called a \emph{chronological list of forces}.

It was shown in \cite{MZC2014} that controllability of $(\mathcal{G};V_{L})$ can be characterized in terms of zero forcing sets.
\begin{prp} \label{p:MCT}
	Let $\mathcal{G} = (V,E)$ be a  directed graph and $V_{L} \subseteq V$ be the leader set. Then, $(\mathcal{G};V_{L})$ is controllable 
	 if and only if $V_{L}$ is a zero forcing set.
\end{prp}

\subsection{Balancing Set and Controllability of $(\mathcal{G}(W);V_{L})$}
Consider the  weighted  graph $\mathcal{G}(W) = (V,E,W)$ associated with the  directed graph $\mathcal{G} = (V,E)$ and weighted adjacency matrix $W \in \mathcal{W}(\mathcal{G})$. For $i = 1, \ldots, n$, let $x_{i}$ be a variable assigned to vertex $i$. 
Assume that for a given subset of vertices $C \subseteq V$, $x_{j} = 0$ for all $ j \in C$.  We call $C$ {\em the set of zero vertices}. 
The  values of the other vertices of $\mathcal{G}(W)$ are initially undetermined.
To every vertex $j \in C$, we assign a so called {\em balance equation}:
\begin{equation}\label{e:BE}
	\sum_{k \in N_{V \setminus C}(\{j\})} x_{k}W_{kj} = 0.
\end{equation}
Note that for weighted  undirected graphs, in which case $W = W^{T}$, the balance equation \eqref{e:BE} coincides with the one introduced in \cite{MHM2018}.
Assume that there is a subset of zero vertices $X \subseteq C$ such that the  system of $|X|$ balance equations corresponding to the vertices in $X$ implies that $x_{k} = 0$ for all $k \in Y$ with $C \cap Y = \emptyset$. 
The updated set of zero vertices is now defined as $C' = C \cup Y$. In  this case, we say that {\em zeros extend from $X$ to  $Y$},  written as $X \xrightarrow{z} Y$.

This one step procedure of making the values of possibly additional vertices equal to zero is called {\em the zero extension rule\/}.

Define  the {\em derived set\/} $\mathcal{D}_{z}(C)$ to be the set of zero vertices obtained after repeated application of the zero extension rule until no more zero vertices can be added. Although not explicitly stated in \cite{MHM2018}, it can be shown that the derived set is uniquely defined, in the sense that it does not depend on the particular zero extensions that are applied to the original set of zero vertices $C$.
 An initial zero vertex set $C \subseteq V$ is called a {\em balancing set} if the derived set $\mathcal{D}_{z}(C)$ is $V$. Given a balancing set, one can list the zero extensions in the order in which they were performed.  Such a list is called a \emph{chronological list of zero extensions.}
 
A necessary and sufficient condition for strong structural controllability with respect to $\mathcal{Q}_{W}(\mathcal{G})$  for the special case that $W = W^T$ was given in \cite{MHM2018}:
\begin{prp} \label{BSC1}
	Let $\mathcal{G}$ be a simple undirected graph, $V_{L} \subseteq V$ be the leader set and  $W \in \mathcal{W}(\mathcal{G})$ be a weighted adjacency matrix with $W = W^{T}$.
	Then $(\mathcal{G}(W);V_{L})$ is controllable if and only if $V_{L}$ is a balancing set. 
\end{prp}

\section{Problem formulation} 

In this section we will introduce the main problem that is considered in this paper. 
At the end of the section, we will also formulate two preliminary results that will be needed in the sequel. 
In order to proceed, we will now first formalize that the weights of a priori given edges in the network graph are constrained to be equal. This can be expressed as a condition that some of the off-diagonal entries in the matrices belonging to the qualitative class $\mathcal{Q}(\mathcal{G})$  are  equal. 
To do this, we introduce a partition
\[
	\pi = \{E_{1},E_{2},\ldots,E_{k}\}
\]
of the edge set $E$ into disjoint subsets $E_r$ whose union is the entire edge set $E$. 
The edges in a given cell $E_r$ are constrained to have identical weights. 
We then define the {\em colored qualitative class} associated with $\pi$ by
\[\begin{split}
\mathcal{Q}_{\pi}(\mathcal{G}) = &\{A \in \mathcal{Q}(\mathcal{G}) \mid  A_{ij} = A_{kl}  \\ &  \mbox{ if } (j,i), (l,k) \in E_{r} \mbox{ for some } r\}.
\end{split} \]

In order to visualize the partition $\pi$ of the edge set in the graph, two edges in the same cell $E_{r}$ are said to have the same color. 
The colors will be denoted by the symbols $c_1, c_2, \ldots, c_k$ and the edges in cell $E_r$ are said to have color $c_r$. This leads to the notion of  {\em colored graph}. 
A {\em colored  graph} is a directed graph together with a partition $\pi$ of the edge set, which is denoted by $\mathcal{G}(\pi) = (V,E,\pi)$.  

In the sequel, sometimes the symbols $c_i$ will also be used to denote independent nonzero variables. 
A set of real values obtained by assigning to each of these variables $c_i$ a particular real value is called a \emph{realization} of the color set.

\begin{ex} \label{ex:coloredgraph}
	Consider the colored graph $\mathcal{G}(\pi) = (V,E,\pi)$ associated  with the directed graph $\mathcal{G} = (V,E)$ and edge partition $\pi =\{E_{1},E_{2},E_{3}\}$, where $ E_{1}= \{(1,4), (1,6)\}$, $E_{2}=\{(2,4),(2,5)\}$  and  $E_{3} = \{(3,5), (3,6)\}$ as depicted in Figure \ref{g:SDCG}.
		\begin{figure}[h!]
		\centering
		\begin{tikzpicture}[scale=0.5]
		\tikzset{VertexStyle1/.style = {shape = circle,
				ball color = white!100!black,
				text = black,
				inner sep = 2pt,
				outer sep = 0pt,
				minimum size = 10 pt},
			edge/.style={->,> = latex', text = black}
		}
		\tikzset{VertexStyle2/.style = {shape = circle,
				ball color = black!80!yellow,
				text = white,
				inner sep = 2pt,
				outer sep = 0pt,
				minimum size = 10 pt}}
		\node[VertexStyle2](1) at (0,3) {$1$};
		\node[VertexStyle2](2) at (-2,0) {$2$};
		\node[VertexStyle2](3) at (2,0) {$3$};
		\node[VertexStyle1](4) at (-5,-3) {$4$};
		\node[VertexStyle1](5) at (0,-3) {$5$};
		\node[VertexStyle1](6) at (5,-3) {$6$};
		\Edge[ style = {->,> = latex'},color=green, label = $c_{2}$,labelstyle={inner sep=0pt}](2)(4);
		\Edge[  style = {->,> = latex',pos = 0.5},color=green, label = $c_{2}$,labelstyle={inner sep=0pt}](2)(5);
		\Edge[style = {->,> = latex',pos = 0.5},color= blue, label = $c_{3}$, labelstyle={inner sep=0pt}](3)(5);
		\Edge[ style = {->,> = latex'},color= blue , label = $c_{3}$,labelstyle={inner sep=0pt}](3)(6);
		\Edge[style={bend left,->,> = latex'},color = red, label = $c_{1}$,labelstyle={inner sep=0pt}](1)(6)
		\Edge[style={bend right,->,> = latex'},color = red, label = $c_{1}$,labelstyle={inner sep=0pt}](1)(4)
		\Edge[style={bend right,->,> = latex'},color = red, label = $c_{1}$,labelstyle={inner sep=0pt}](4)(5)
		\Edge[style={bend left,->,> = latex'},color = blue, label = $c_{3}$,labelstyle={inner sep=0pt}](6)(5)
		\end{tikzpicture}
		\caption{A colored  directed graph with leader set $\{1,2,3\}$.}
		\label{g:SDCG}
	\end{figure}
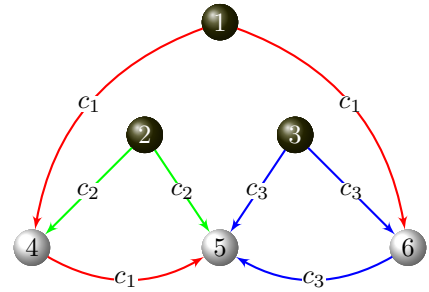
Edges having the same color means that the weight of these edges are constrained to be equal. 
In this example, the edges in $E_1$ have color $c_1$ (blue), those in $E_2$ have color $c_2$ (green), and those in $E_3$ have color $c_3$ (red). 
The corresponding colored qualitative class consists of all matrices of the form
\[\begin{bmatrix}\lambda_1& 0 & 0 & 0 & 0 & 0\\
0 & \lambda_2& 0 & 0 & 0 & 0\\
0 & 0&  \lambda_3& 0 & 0 & 0 \\
c_1 & c_2 & 0 & \lambda_4& 0 & 0 \\
0 & c_2 & c_3 & c_1 &\lambda_5& c_3 \\
c_1&0 & c_3 &0& 0  & \lambda_6
 \end{bmatrix}\]
where $\lambda_i$ is an arbitrary real number for $i=1,2,\ldots,6$ and $c_i$ is an arbitrary {\em nonzero} 
real number for $i=1,2,3.$
\end{ex}
\smallskip

Given a colored directed graph $\mathcal{G}(\pi) = (V,E,\pi)$ with edge partition $\pi = \{E_{1}, E_{2}, \ldots, E_{k}\}$, we define the corresponding family of weighted adjacency matrices
\[\begin{split}
\mathcal{W}_{\pi}(\mathcal{G})  := &\{ W \in \mathcal{W}(\mathcal{G}) \mid W_{ij} = W_{kl} \\ &  \mbox{ if } (j,i),(l,k) \in E_{r} \mbox{ for some } r\}.
\end{split} \]
Note that any weighted adjacency matrix $W \in \mathcal{W}_{\pi}(\mathcal{G})$ is associated with a unique realization of the color set. 
Obviously, the colored qualitative class $\mathcal{Q}_{\pi}(\mathcal{G})$ is equal to the union of all the subclasses $\mathcal{Q}_{W}(\mathcal{G})$ with $W \in \mathcal{W}_{\pi}(\mathcal{G})$, i.e, 
\begin{equation}  \label{e:union}
\mathcal{Q}_{\pi}(\mathcal{G}) = \bigcup_{W \in  \mathcal{W}_{\pi}(\mathcal{G})} \mathcal{Q}_{W}(\mathcal{G}).
\end{equation}  

If $(\mathcal{G};V_{L})$ is controllable with respect to $\mathcal{Q}' = \mathcal{Q}_{\pi}(\mathcal{G})$ (see Definition \ref{d:cos}) we will simply say that {\em $(\mathcal{G}(\pi);V_{L})$ is controllable}. 
In that case, we call the system {\em colored strongly structurally controllable}.
For example, the system with graph depicted in Figure \ref{g:SDCG} is colored strongly structurally controllable as will be shown later. 

The aim of this paper is to establish graph-theoretic tests for colored strong structural controllability of a given graph. 
In order to obtain such conditions, we now first make the observation that conditions for strong structural controllability can be expressed in terms of balancing sets. 
Generalizing Proposition  \ref{BSC1} to the case of weighted directed graphs, we have the following lemma:  
\begin{lem} \label{l:BSC}
	Let $\mathcal{G} = (V,E)$ be a directed graph with leader set $V_{L}$ and let $W \in \mathcal{W}(\mathcal{G})$. 
	Then $(\mathcal{G}(W); V_{L})$ is controllable if and only if $V_{L}$ is a balancing set.	
\end{lem}

The proof can be found in the Appendix.

The following lemma follows immediately from Lemma \ref{l:BSC} by noting that \eqref{e:union} holds.
\begin{lem}\label{l:CSSCBS}
Let $\mathcal{G} = (V,E)$ be a directed graph with leader set $V_{L}$ and let $\pi$ be a partition of the edge set. 
Then $(\mathcal{G}(\pi);V_{L})$ is controllable if and only if  $V_{L}$ is a balancing set for all weighted graphs $\mathcal{G}(W) = (V,E,W)$ with $W \in \mathcal{W}_{\pi}(\mathcal{G})$.	
\end{lem}

Obviously, the necessary and sufficient conditions presented in Lemma~\ref{l:CSSCBS} cannot be verified easily, as the set $\mathcal{W}_{\pi}(\mathcal{G})$ contains infinitely many elements. Therefore, we aim at establishing graph-theoretic conditions under which $(\mathcal{G}(\pi);V_{L})$ is controllable. 


\section{Zero Forcing Sets For Colored Graphs}

In order to provide a graph-theoretic condition for colored strong structural controllability, in this section we introduce a new color change rule and then define the corresponding notion of zero forcing set. 
To do this, we first consider colored bipartite graphs and establish a necessary and sufficient graph-theoretic condition for nonsingularity of the associated pattern class.

\subsection{Colored Bipartite Graphs}
Consider the bipartite graph $G = (X,Y,E_{XY})$, where the vertex sets $X$ and $Y$ are given by $X = \{x_{1},x_{2}, \ldots, x_{s}\}$ and $Y = \{y_{1},y_{2}, \ldots, y_{t}\}$. 
We will now introduce the notion of {\em colored} bipartite graph. Let $\pi_{XY} = \{E_{XY}^{1}, E_{XY}^{2}, \ldots, E_{XY}^{\ell}\}$ be a partition of the edge set $E_{XY}$ with associated colors $c_{1}, c_{2}, \ldots, c_{\ell}$.
This partition is now used to formalize that certain entries in the pattern class $\mathcal{P}(G)$ are constrained to be equal. 
Again, the edges in a given cell $E^{r}_{XY}$ are said to have the same color.
The pattern class of  the colored bipartite graph $G(\pi) = (X,Y,E_{XY},\pi_{XY})$ is  then defined as the following set of complex $t \times s$ matrices
\[
\begin{split}
& \mathcal{P}_{\pi}(G) =
\big\{M \in \mathcal{P}(G) \mid M_{ji} = M_{hg} \\ & \mbox{ if }  \{x_{i},y_{j}\},  \{x_{g},y_{h} \}\in  E^{r}_{XY} \mbox{ for some } r\big\}.
\end{split}
\]
Assume now that $|X| = |Y|$ and let $t=|X|$.
Suppose that $p$ is a perfect matching of $G(\pi)$. 
The  {\em spectrum} of $p$ is  defined to be the set of colors (counting multiplicity) of the edges in $p$.
More specifically, if the perfect matching $p$ is given by $p = \big\{\{x_{1},y_{\gamma(1)}\}, \ldots, \{x_{t},y_{\gamma(t)}\}\big\}$, where $\gamma$ denotes a permutation of $(1,2,\ldots,t)$, and $c_{i_1}, c_{i_2}, \ldots, c_{i_t}$ are the respective colors of the edges in $p$, then the spectrum of $p$ is $\{c_{i_1}, c_{i_2}, \ldots, c_{i_t}\}$ where the same color can appear multiple times.
	
In addition, we define the {\em sign} of the perfect matching $p$ as $\sign(p) = (-1)^{m}$, where $m$ is the number of swaps needed to obtain $(\gamma(1), \gamma(2), \ldots , \gamma(t))$ from $(1,2, \ldots, t)$.  
Since every perfect matching  is associated with a unique permutation, with a slight abuse of notation, we sometimes use the perfect matching $p$ to represent its corresponding permutation.

 Two perfect matchings are called  {\em equivalent}  if they have the same spectrum.  
%
%
Obviously this yields a partition of the set of all perfect matchings of $G(\pi)$ into   {\em equivalence classes} of perfect matchings.
We denote these equivalence classes of perfect matchings by $\mathbb{P}_{1},\mathbb{P}_{2},\ldots,\mathbb{P}_{l}$, where perfect matchings in the same class $\mathbb{P}_{i}$  are equivalent. 
Clearly, $\mathbb{P}_{i} \cap \mathbb{P}_{j} = \emptyset$ for $i \neq j$.
Correspondingly, we then define the {\em spectrum of the equivalence class $\mathbb{P}_{i}$} to be the (common) spectrum of the perfect matchings in this class, and denote it by $\spec(\mathbb{P}_{i})$.
Finally, we define the {\em the signature of the equivalence class $\mathbb{P}_{i}$} to be the sum of the signs of all perfect matchings in this class, which is given by
\[\sgn(\mathbb{P}_{i}) = \sum_{p \in \mathbb{P}_{i}} \sign(p).\]
\begin{ex} \label{ex:coloredbg}
	Consider the colored bipartite graph $G(\pi)$ depicted in Figure \ref{g:CBG1}.  
	It contains three perfect matchings, $p_{1}$, $p_{2}$ and $p_{3}$, respectively, depicted in Figure \ref{g:CBG}(b)-(d). 
	Clearly, $p_{1}$ and $p_{3}$ are equivalent. 
	The equivalence classes of perfect matchings are then $\mathbb{P}_{1} = \{p_{1},p_{3}\}$ and $\mathbb{P}_{2} = \{p_{2}\}$. 
	Clearly, $\sgn(\mathbb{P}_{1}) = 0$ and $\sgn(\mathbb{P}_{2}) = -1$.
\end{ex}
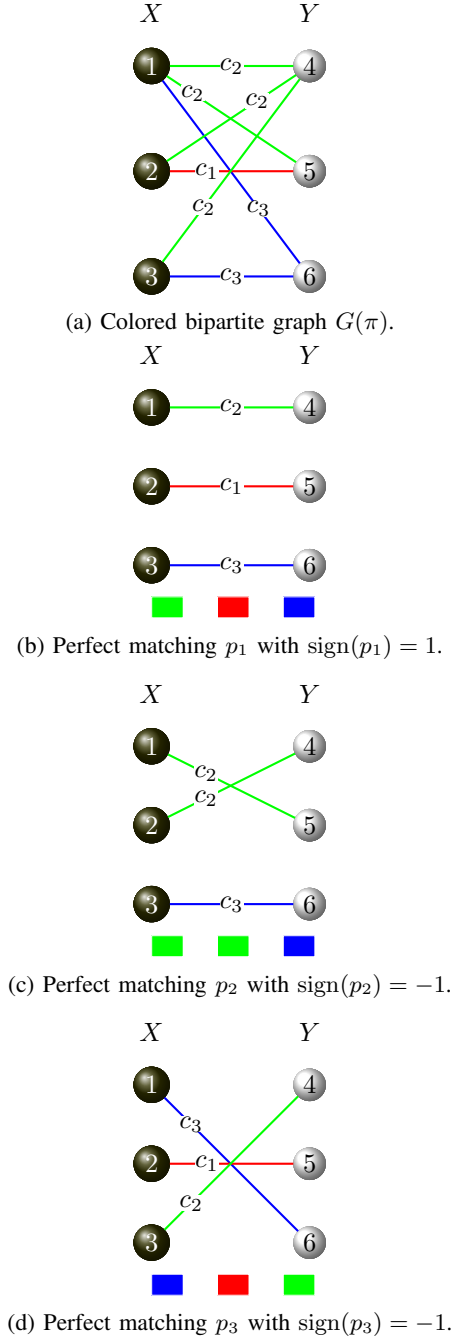
\begin{figure}[h!]
	\centering
	\begin{subfigure}{0.4\textwidth}\label{f:1}
		\centering
		\begin{tikzpicture}[scale=0.7]
		\tikzset{VertexStyle1/.style = {shape = circle,
				ball color = white!100!black,
				text = black,
				inner sep = 1.5pt,
				outer sep = 0pt,
				minimum size = 6 pt},
			edge/.style={->,> = latex', text = black}
		}
		\tikzset{VertexStyle2/.style = {shape = circle,
				ball color = black!80!yellow,
				text = white,
				inner sep = 2pt,
				outer sep = 0pt,
				minimum size = 10 pt}}
		\node[VertexStyle2](1) at (-1.5,2) {$1$};
		\node[VertexStyle2](2) at (-1.5,0) {$2$};
		\node[VertexStyle2](3) at (-1.5,-2) {$3$};
		\node[VertexStyle1](4) at (1.5,2) {$4$};
		\node[VertexStyle1](5) at (1.5,0) {$5$};
		\node[VertexStyle1](6) at (1.5,-2) {$6$};
		\node[](7) at (-1.5,3){$X$};
		\node[](8) at (1.5,3){$Y$};
		\Edge[style = {,> = latex',pos = 0.5},color=green
		, label = $c_{2}$,labelstyle={inner sep=0pt}](1)(4);
		\Edge[ style = {,> = latex',pos = 0.2},color=green
		, label = $c_{2}$,labelstyle={inner sep=0pt}](1)(5);
		\Edge[style = {,> = latex',pos = 0.7},color=blue
		, label = $c_{3}$,labelstyle={inner sep=0pt}](1)(6);
		\Edge[ style = {,> = latex',pos = 0.7},color=green
		, label = $c_{2}$,labelstyle={inner sep=0pt}](2)(4);
		\Edge[ style = {,> = latex',pos = 0.3},color=red
		, label = $c_{1}$,labelstyle={inner sep=0pt}](2)(5);
		\Edge[style = {,> = latex',pos = 0.3},color=green
		, label = $c_{2}$,labelstyle={inner sep=0pt}](3)(4);
		\Edge[style = {,> = latex',pos = 0.5},color=blue
		, label = $c_{3}$,labelstyle={inner sep=0pt}](3)(6);
		\end{tikzpicture}
		\caption{Colored bipartite graph $G(\pi)$.}
		\label{g:CBG1}
	\end{subfigure}
	\vspace{.2cm}
	\begin{subfigure}{0.4\textwidth}\label{f:2}
		\centering
		\begin{tikzpicture}[scale=0.7]
		\tikzset{VertexStyle1/.style = {shape = circle,
				ball color = white!100!black,
				text = black,
				inner sep = 1.5pt,
				outer sep = 0pt,
				minimum size = 8 pt},
			edge/.style={->,> = latex', text = black}
		}
		\tikzset{VertexStyle2/.style = {shape = circle,
				ball color = black!80!yellow,
				text = white,
				inner sep = 2pt,
				outer sep = 0pt,
				minimum size = 10 pt}}
		\node[VertexStyle2](1) at (-1.5,2) {$1$};
		\node[VertexStyle2](2) at (-1.5,0.5) {$2$};
		\node[VertexStyle2](3) at (-1.5,-1) {$3$};
		\node[VertexStyle1](4) at (1.5,2) {$4$};
		\node[VertexStyle1](5) at (1.5,0.5) {$5$};
		\node[VertexStyle1](6) at (1.5,-1) {$6$};
		\node[](7) at (-1.5,3){{$X$}};
		\node[](8) at (1.5,3){{$Y$}};
		\Edge[ style = {,> = latex',pos = 0.5},color=green
		, label = $c_{2}$,labelstyle={inner sep=0pt}](1)(4);
		\Edge[ style = {,> = latex',pos = 0.5},color=red
		, label = $c_{1}$,labelstyle={inner sep=0pt}](2)(5);
		\Edge[style = {,> = latex',pos = 0.5},color=blue
		, label = $c_{3}$,labelstyle={inner sep=0pt}](3)(6);
		\fill[fill=green,draw=white] (-1.5,-2) rectangle +(6mm,4mm);
		\fill[fill=red,draw=white] (-0.25,-2) rectangle +(6mm,4mm);
		\fill[fill=blue,draw=white] (1.0,-2) rectangle +(6mm,4mm);
		\end{tikzpicture}
		\caption{Perfect matching $p_{1}$ with $\sign(p_{1}) = 1$.}
		\label{g:CBG2}
	\end{subfigure}
	\vspace{.2cm}
	\begin{subfigure}{0.4\textwidth}\label{f:3}
		\centering
		\begin{tikzpicture}[scale=0.7]
		\tikzset{VertexStyle1/.style = {shape = circle,
				ball color = white!100!black,
				text = black,
				inner sep = 1.5pt,
				outer sep = 0pt,
				minimum size = 8 pt},
			edge/.style={->,> = latex', text = black}
		}
		\tikzset{VertexStyle2/.style = {shape = circle,
				ball color = black!80!yellow,
				text = white,
				inner sep = 2pt,
				outer sep = 0pt,
				minimum size = 10 pt}}
		\node[VertexStyle2](1) at (-1.5,2) {$1$};
		\node[VertexStyle2](2) at (-1.5,0.5) {$2$};
		\node[VertexStyle2](3) at (-1.5,-1) {$3$};
		\node[VertexStyle1](4) at (1.5,2) {$4$};
		\node[VertexStyle1](5) at (1.5,0.5) {$5$};
		\node[VertexStyle1](6) at (1.5,-1) {$6$};
		\node[](7) at (-1.5,3){{$X$}};
		\node[](8) at (1.5,3){{$Y$}};
		\Edge[style = {,> = latex',pos = 0.3},color=green
		, label = $c_{2}$,labelstyle={inner sep=0pt}](1)(5);
		\Edge[ style = {,> = latex',pos = 0.3
		},color=green
		, label = $c_{2}$,labelstyle={inner sep=0pt}](2)(4);
		\Edge[style = {,> = latex',pos = 0.5},color=blue
		, label = $c_{3}$,labelstyle={inner sep=0pt}](3)(6);
		\fill[fill=green,draw=white] (-1.5,-2) rectangle +(6mm,4mm);
		\fill[fill=green,draw=white] (-0.25,-2) rectangle +(6mm,4mm);
		\fill[fill=blue,draw=white] (1.0,-2) rectangle +(6mm,4mm);
		\end{tikzpicture}
		\caption{Perfect matching $p_{2}$ with $\sign(p_{2}) = -1$.}
		\label{g:CBG3}
	\end{subfigure}
	\vspace{.2cm}
	\begin{subfigure}{0.4\textwidth}\label{f:4}
		\centering
		\begin{tikzpicture}[scale=0.7]
		\tikzset{VertexStyle1/.style = {shape = circle,
				ball color = white!100!black,
				text = black,
				inner sep = 1.5pt,
				outer sep = 0pt,
				minimum size = 8 pt},
			edge/.style={->,> = latex', text = black}
		}
		\tikzset{VertexStyle2/.style = {shape = circle,
				ball color = black!80!yellow,
				text = white,
				inner sep = 2pt,
				outer sep = 0pt,
				minimum size = 10 pt}}
		\node[VertexStyle2](1) at (-1.5,2) {$1$};
		\node[VertexStyle2](2) at (-1.5,0.5) {$2$};
		\node[VertexStyle2](3) at (-1.5,-1) {$3$};
		\node[VertexStyle1](4) at (1.5,2) {$4$};
		\node[VertexStyle1](5) at (1.5,0.5) {$5$};
		\node[VertexStyle1](6) at (1.5,-1) {$6$};
		\node[](7) at (-1.5,3){$X$};
		\node[](8) at (1.5,3){$Y$};
		\Edge[ style = {,> = latex',pos = 0.2},color=blue
		, label = $c_{3}$,labelstyle={inner sep=0pt}](1)(6);
		\Edge[ style = {,> = latex',pos = 0.3},color=red
		, label = $c_{1}$,labelstyle={inner sep=0pt}](2)(5);
		\Edge[style = {,> = latex',pos = 0.2},color=green
		, label = $c_{2}$,labelstyle={inner sep=0pt}](3)(4);
		\fill[fill=blue,draw=white] (-1.5,-2) rectangle +(6mm,4mm);
		\fill[fill=red,draw=white] (-0.25,-2) rectangle +(6mm,4mm);
		\fill[fill=green,draw=white] (1.0,-2) rectangle +(6mm,4mm);
		\end{tikzpicture}
		\caption{Perfect matching $p_{3}$ with $\sign(p_{3}) = -1$.}
		\label{g:CBG4}
	\end{subfigure}
	\caption{Example of a colored bipartite graph and its perfect matchings.}
	\label{g:CBG}
\end{figure}

We are now ready to state a necessary and sufficient condition for nonsingularity of all matrices in the colored pattern class $\mathcal{P}_{\pi}(G)$.

\begin{thm}\label{t:NoP}
	Let $G(\pi) = (X,Y,E_{XY},\pi_{XY})$  be a colored bipartite graph and $|X| = |Y|$.  
	Then, all matrices in $\mathcal{P}_{\pi}(G)$ are nonsingular if and only if there exists at least one perfect matching and exactly one equivalence class of perfect matchings has nonzero signature.
\end{thm}
\begin{IEEEproof}
Denote the cardinality of $X$ and $Y$ by $t$.
Let $A \in \mathcal{P}_{\pi}(G)$. 
By the Leibniz Formula for the determinant, we have \[\det(A) = \sum_{\gamma} \sign(\gamma) \prod_{i = 1}^{t}A_{i \gamma(i)},\]
where the sum ranges over all permutations $\gamma$ of $(1,2,\ldots,t)$ and where 
$\sign(\gamma) = (-1)^m$ with $m$ the number of swaps needed to obtain $(\gamma(1), \gamma(2), \ldots , \gamma(t))$ from $(1,2, \ldots, t)$.  
Note that $\prod_{i = 1}^{t}A_{i \gamma(i)} \neq 0$ if and only if there exists at least one perfect matching $ p =\{\{x_{1},y_{\gamma(1)}\}, \ldots, \{x_{|X|},y_{\gamma(t)}\}\}$ in $G(\pi)$. 
In that case, we have
\[
\det(A) = \sum_{p} \sign(p) \prod_{i = 1}^{t}A_{i p(i)} ,
\]
where $p$ ranges over all perfect matchings and $\sign(p)$ denotes the sign of the perfect matching (we now identify perfect matchings with their permutations).
Suppose now there are $l$ equivalence classes of perfect matchings $\mathbb{P}_{1},\mathbb{P}_{2},\ldots,\mathbb{P}_{l}$. 
Then we obtain
\begin{equation} \label{e:det}
  \det(A) = \sum_{j=1}^l \left( \sgn(\mathbb{P}_{j}) \prod_{i=1}^t A_{i p(i)} \right),
\end{equation}
where, for $j = 1,2, \ldots l$, in the product appearing in the $j$th term, $p$ is an arbitrary matching in $\mathbb{P}_j$.
We will now prove the `if' part. 
Assume that there exists at least one perfect matching, and exactly one equivalence class of perfect matchings has nonzero signature. 
Without loss of generality,  assume that the equivalence class $\mathbb{P}_{1}$ has nonzero signature. 
Obviously, for every $A \in \mathcal{P}_{\pi}(G)$, we then have
\[
	\det(A) =  \sgn(\mathbb{P}_{1}) \prod_{i=1}^t A_{i p(i)} \neq 0 ,
\]
where $p \in \mathbb{P}_1$ is arbitrary, in other words, every $A \in \mathcal{P}_{\pi}(G)$ is nonsingular.
	
Next, we prove the `only if' part. 
For this, assume that all $A \in \mathcal{P}_{\pi}(G)$ are nonsingular, but one of the following holds:
\begin{itemize}
\item[(i)] there does not exist any perfect matching,
\item[(ii)] no equivalence class of perfect matchings with nonzero signature exists, or
\item[(iii)] there exist at least two equivalence classes of perfect matchings with nonzero signature.
\end{itemize}
We will show that all these cases lead to contradiction.

In case (i), we must obviously have $\det(A) = 0$ for any $A \in \mathcal{P}_{\pi}(G)$ which gives a contradiction. 
For case (ii), it follows from \eqref{e:det} that $\det(A) = 0$ since all equivalence classes have zero signature.  
Therefore, we reach a contradiction again. 
Finally, consider case (iii). 
Without loss of generality, assume $\mathbb{P}_{1}$ and $\mathbb{P}_{2}$ have nonzero signature. 
The signatures of the remaining equivalence classes can be either zero or nonzero.  
In the sequel we associate the colors $c_1, c_2, \ldots, c_{\ell}$ of the cells $E_{XY}^1, E_{XY}^2, \ldots E_{XY}^{\ell}$ with independent, nonzero, variables $c_1, c_2, \ldots, c_{\ell}$ that can take values in $\mathbb{C}$. 
The spectrum of an equivalence class $\mathbb{P}_{j}$ then uniquely determines a monomial $c_1^{i_1}c_2^{i_2} \ldots c_{\ell}^{i_{\ell}}$, where the powers $i_1, i_2, \ldots i_{k}$ corres\-pond to the multiplicities of the colors $c_1, c_2, \ldots ,c_{\ell}$ in the perfect matchings in $\mathbb{P}_{j}$. 
We also identify each entry of a matrix $A$ in $\mathcal{P}_{\pi}(G)$ with the color of its corresponding edge. 
In particular, for such $A$ we have
\[A_{ij} = 
	\left\{ 
		\begin{split}
		&  c_{r} & \mbox{  if } (j,i) \in E_{r} \mbox{ for some } r,\\
		&  0      & \mbox{  otherwise,}
		\end{split}
		\right.
\]
From the expression \eqref{e:det} for the determinant of $A$ it can be seen that the perfect matchings in the equivalence class $\mathbb{P}_{j}$ yield a contribution $\sgn(\mathbb{P}_{j}) c_1^{i_1}c_2^{i_2} \ldots c_{\ell}^{i_\ell}$, where the degrees correspond to the multiplicities of the colors of the perfect matchings in $\mathbb{P}_{j}$. 
By assumption we have that $\spec(\mathbb{P}_{1})$ and $\spec(\mathbb{P}_{2})$ are not equal. 
Without loss of generality, we assume that the multiplicity of $c_{1}$ as an element of $\spec(\mathbb{P}_{1})$ is unequal to the multiplicity of $c_{1}$ as an element of $\spec(\mathbb{P}_{2})$. 
Denote these multiplicities by $j_{1}$ and $j_{2}$, respectively,  with $j_1 \neq j_2$. 
Then for all values of $c_{2}, \ldots, c_\ell$, the determinant of $A$ has the form
\begin{equation} \label{e:det1}
\det(A) = \sgn(\mathbb{P}_{1})a_{1}c_{1}^{j_{1}}+\sgn(\mathbb{P}_{2})a_{2}c_{1}^{j_{2}}+f(c_{1}),
\end{equation}
where $a_{1}$ and $a_{2}$ depend on $c_{2}, \ldots, c_k$ and $f(c_{1})$ is a polynomial in $c_{1}$. 
The polynomial $f(c_1)$ corresponds to the remaining equivalence classes. 
It can happen that some of these equivalence classes also contain the color $c_1$ in their spectrum with multiplicity $j_1$ or $j_2$.  
By moving the corresponding monomials to the first two terms in \eqref{e:det1} we obtain
\begin{equation} \label{e:det2}
\det(A) = b_{1}c_{1}^{j_{1}} + b_{2}c_{1}^{j_{2}} + f'(c_{1}),
\end{equation}
with  $b_1$ and $b_2$ depending on $c_{2}, \ldots, c_k$. 
Note that the first term in \eqref{e:det2} corresponds to the equivalence classes containing $c_1$ in their spectrum with multiplicity $j_1$, and likewise the second term with multiplicity $j_2$. 
The remaining polynomial $f'(c_1)$ does not contain monomials with $c_{1}^{j_{1}}$ and $c_{1}^{j_{2}}$. 
It is now easily verified  that nonzero $c_{2}, \ldots, c_\ell$ can be chosen such that $b_1 \neq 0$ and $b_2 \neq 0$.  
By the fundamental theorem of algebra we then have that the polynomial equation $b_{1}c_{1}^{j_{1}}+  b_{2}c_{1}^{j_{2}} + f'(c_{1})=0$ has at least one nonzero root, since both $b_{1}$ and  $b_{2}$ are nonzero. 
This implies that for some choice of nonzero complex values $c_{1}, c_{2}, \ldots, c_\ell$ we have $\det(A) =0$. 
In other words, not all $A \in \mathcal{P}_{\pi}(G)$ are nonsingular. 
This is a contradiction.  
\end{IEEEproof}
\begin{ex}
For the colored bipartite graph in Figure \ref{g:CBG1}, the pattern class consists of 
all matrices of the form
\[\begin{bmatrix}
c_2 & c_2 & c_2 \\
c_2 & c_1 & 0    \\
c_3 & 0   &  c_3 \\
\end{bmatrix}\]
where $c_1,c_2$ and $c_3$ are arbitrary nonzero complex numbers. In Example \ref{ex:coloredbg}
we saw that there is exactly one equivalence class of perfect matchings with nonzero signature. By Theorem 8 we thus conclude that all these matrices are nonsingular.
\end{ex}

\subsection{Color Change Rule and Zero Forcing Sets}

In this subsection, we will introduce a tailor-made zero forcing notion for colored graphs.
Let $\mathcal{G}(\pi) = (V,E,\pi)$ be a colored directed graph with $\pi = \{E_{1}, E_{2}, \ldots, E_{k}\}$ the partition of $E$.
For given disjoint subsets $X = \{x_{1},x_{2}, \ldots, x_{s}\}$ and $Y = \{y_{1},y_{2}, \ldots, y_{t}\}$ of $V$, we define an associated colored bipartite graph $G(\pi) = (X,Y,E_{XY}, \pi_{XY})$ as follows:
\[E_{XY} := \{ \{x_i,y_j\} \mid (x_i,y_j) \in E, ~x_i \in X,~ y_j \in Y \}.\]
Obviously, the partition $\pi$ induces a partition  $\pi_{XY}$ of $E_{XY}$ by defining
\[E_{XY}^{r} := \{ \{x_i,y_j\} \in E_{XY} \mid (x_i,y_j) \in E_r \}, ~r =1,2 \ldots, k.\]
Note that for some $r$, this set might be empty. Removing these, we get a partition
\[\pi_{XY} = \{E_{XY}^{i_{1}}, E_{XY}^{i_{1}}, \ldots, E_{XY}^{i_{\ell}}\}\]
of $E_{XY}$, with associated colors $c_{i_{1}}, c_{i_{2}}, \ldots, c_{i_{\ell}}$, with $\ell \leq k$.
Without loss of generality we renumber $c_{i_{1}}, c_{i_{2}}, \ldots, c_{i_{\ell}}$ as $c_{1}, c_{2}, \ldots, c_{\ell}$  and the edges in cell $E^{r}_{XY}$ are said to have color $c_r$.

As before, a subset $C$ of $V$ is called a {\em coloring set} if the vertices in $C$ are initially colored black and those in $V \setminus C$ initially colored white. We will now define the notion of {\em color-perfect white neighbor}.

\begin{defn}\label{d:CPN}
Let $X \subseteq C$ and $Y \subseteq V$ with $|Y| = |X|$.  We call $Y$ a {\em color-perfect white neighbor\/} of $X$ if 
\begin{enumerate}
\item
$Y = N_{V \setminus C}(X)$, i.e. $Y$ is equal to the set of white out-neighbors of $X$, and 
\item
in the associated colored bipartite graph $G=(X,Y,E_{XY},\pi_{XY})$ there exists a perfect matching and exactly one equivalence class of perfect matchings has nonzero signature.
\end{enumerate}
\end{defn}

Based on the notion of color-perfect white neighbor, we now introduce the following color change rule:
if $X\subseteq C$ and $Y$ is a color-perfect white neighbor of $X$, then we change the color of all vertices in $Y$ to black, and write $X {\xrightarrow{c}} Y$. Such a color change is called a {\em force}. We define {\em a derived set\/} $\mathcal{D}_c(C)$ as a set of black vertices obtained after repeated application of the color change rule, until no more changes are possible. In contrast with the original color change rule (see Section \ref{s:ZFS}), under our new color change rule derived sets will no longer be uniquely defined, and may depend on the particular list of forces that is applied to the original coloring set $C$. This is illustrated by Example \ref{ex:counterexample} in the Appendix.

A coloring set $C \subseteq V$ is  called a {\em zero forcing set for} $\mathcal{G}(\pi)$ if there exits a derived set $\mathcal{D}_c(C)$ such that $\mathcal{D}_c(C) = V$.

Before illustrating the new color change rule, we remark on its relation to the one defined earlier.
\begin{remark}
Given a directed graph $\mathcal{G} = (V,E)$, one can obtain a colored graph $\mathcal{G}(\pi) = (V,E,\pi)$ by assigning to every edge a different color, i.e., $|\pi| = |E|$. Clearly, the colored qualitative class $\mathcal{Q}_{\pi}(\mathcal{G})$ coincides with the qualitative class $ \mathcal{Q}(\mathcal{G})$. In addition, the original color change rule for $\mathcal{G}$ introduced in Section \ref{s:ZFS} can be seen to be a special case of the new one for $\mathcal{G}(\pi)$. This observation in mind, we will use the same terminology for these two color change rules and it will be clear from the context which one is employed. 
\end{remark}
We now illustrate the new color change rule by means of an example.
\begin{ex} \label{ex:ZFS}
	\begin{figure}[h!]
		\centering
		\begin{subfigure}{0.4\textwidth}
			\centering
			\begin{tikzpicture}[scale=0.4]
			\tikzset{VertexStyle1/.style = {shape = circle,
					ball color = white!100!black,
					text = black,
					inner sep = 1.5pt,
					outer sep = 0pt,
					minimum size = 8 pt},
				edge/.style={->,> = latex', text = black}
			}
			\tikzset{VertexStyle2/.style = {shape = circle,
					ball color = black!80!yellow,
					text = white,
					inner sep = 2pt,
					outer sep = 0pt,
					minimum size = 10 pt}}
			\node[VertexStyle2](1) at (-3,4) {$1$};
			\node[VertexStyle2](2) at (-3,0) {$2$};
			\node[VertexStyle2](3) at (-3,-4) {$3$};
			\node[VertexStyle1](4) at (3,4) {$4$};
			\node[VertexStyle1](5) at (3,0) {$5$};
			\node[VertexStyle1](6) at (3,-4) {$6$};
			\node[VertexStyle1](7) at (9,4) {$7$};
			\node[VertexStyle1](8) at (9,0) {$8$};
			\node[VertexStyle1](9) at (9,-4) {$9$};
			\Edge[style = {->,> = latex',pos = 0.5},color=green
			, label = $c_{2}$,labelstyle={inner sep=0pt}](1)(4);
			\Edge[ style = {->,> = latex',pos = 0.2},color=green
			, label = $c_{2}$,labelstyle={inner sep=0pt}](1)(5);
			\Edge[style = {->,> = latex',pos = 0.7},color=blue
			, label = $c_{3}$,labelstyle={inner sep=0pt}](1)(6);
			\Edge[ style = {->,> = latex',pos = 0.7},color=green
			, label = $c_{2}$,labelstyle={inner sep=0pt}](2)(4);
			\Edge[ style = {->,> = latex',pos = 0.3},color=red
			, label = $c_{1}$,labelstyle={inner sep=0pt}](2)(5);
			\Edge[style = {->,> = latex',pos = 0.3},color=green
			, label = $c_{2}$,labelstyle={inner sep=0pt}](3)(4);
			\Edge[style = {->,> = latex',pos = 0.5},color=blue
			, label = $c_{3}$,labelstyle={inner sep=0pt}](3)(6);
			
			\Edge[style = {->,> = latex',pos = 0.5},color=green
			, label = $c_{2}$,labelstyle={inner sep=0pt}](4)(7);
			\Edge[ style = {->,> = latex',pos = 0.2},color=green
			, label = $c_{2}$,labelstyle={inner sep=0pt}](4)(8);
			\Edge[style = {->,> = latex',pos = 0.3},color=blue
			, label = $c_{3}$,labelstyle={inner sep=0pt}](5)(8);
			\Edge[ style = {->,> = latex',pos = 0.7},color=blue
			, label = $c_{3}$,labelstyle={inner sep=0pt}](5)(9);
			\Edge[ style = {->,> = latex',pos = 0.2},color=red
			, label = $c_{1}$,labelstyle={inner sep=0pt}](6)(7);
			\Edge[style = {->,> = latex',pos = 0.3},color=red
			, label = $c_{1}$,labelstyle={inner sep=0pt}](6)(9);
			
			\Edge[style={bend right,->,> = latex'},color = blue, label = $c_{3}$,labelstyle={inner sep=0pt}](1)(3)
			\Edge[style={bend left,->,> = latex'},color = blue, label = $c_{3}$,labelstyle={inner sep=0pt}](8)(9)
			\Edge[style={->,> = latex'},color = red, label = $c_{1}$,labelstyle={inner sep=0pt}](4)(5)
			\end{tikzpicture}
			\caption{Initial.}
			\label{g:CZFS1}
		\end{subfigure}
		\medskip
		
		\begin{subfigure}{0.4\textwidth}
			\centering
			\begin{tikzpicture}[scale=0.4]
			\tikzset{VertexStyle1/.style = {shape = circle,
					ball color = white!100!black,
					text = black,
					inner sep = 1.5pt,
					outer sep = 0pt,
					minimum size = 8 pt},
				edge/.style={->,> = latex', text = black}
			}
			\tikzset{VertexStyle2/.style = {shape = circle,
					ball color = black!80!yellow,
					text = white,
					inner sep = 2pt,
					outer sep = 0pt,
					minimum size = 10 pt}}
			\node[VertexStyle2](1) at (-3,4) {$1$};
			\node[VertexStyle2](2) at (-3,0) {$2$};
			\node[VertexStyle2](3) at (-3,-4) {$3$};
			\node[VertexStyle2](4) at (3,4) {$4$};
			\node[VertexStyle2](5) at (3,0) {$5$};
			\node[VertexStyle2](6) at (3,-4) {$6$};
			\node[VertexStyle1](7) at (9,4) {$7$};
			\node[VertexStyle1](8) at (9,0) {$8$};
			\node[VertexStyle1](9) at (9,-4) {$9$};
			\Edge[style = {->,> = latex',pos = 0.5},color=green
			, label = $c_{2}$,labelstyle={inner sep=0pt}](1)(4);
			\Edge[ style = {->,> = latex',pos = 0.2},color=green
			, label = $c_{2}$,labelstyle={inner sep=0pt}](1)(5);
			\Edge[style = {->,> = latex',pos = 0.7},color=blue
			, label = $c_{3}$,labelstyle={inner sep=0pt}](1)(6);
			\Edge[ style = {->,> = latex',pos = 0.7},color=green
			, label = $c_{2}$,labelstyle={inner sep=0pt}](2)(4);
			\Edge[ style = {->,> = latex',pos = 0.3},color=red
			, label = $c_{1}$,labelstyle={inner sep=0pt}](2)(5);
			\Edge[style = {->,> = latex',pos = 0.3},color=green
			, label = $c_{2}$,labelstyle={inner sep=0pt}](3)(4);
			\Edge[style = {->,> = latex',pos = 0.5},color=blue
			, label = $c_{3}$,labelstyle={inner sep=0pt}](3)(6);
			
			\Edge[style = {->,> = latex',pos = 0.5},color=green
			, label = $c_{2}$,labelstyle={inner sep=0pt}](4)(7);
			\Edge[ style = {->,> = latex',pos = 0.2},color=green
			, label = $c_{2}$,labelstyle={inner sep=0pt}](4)(8);
			\Edge[style = {->,> = latex',pos = 0.3},color=blue
			, label = $c_{3}$,labelstyle={inner sep=0pt}](5)(8);
			\Edge[ style = {->,> = latex',pos = 0.7},color=blue
			, label = $c_{3}$,labelstyle={inner sep=0pt}](5)(9);
			\Edge[ style = {->,> = latex',pos = 0.2},color=red
			, label = $c_{1}$,labelstyle={inner sep=0pt}](6)(7);
			\Edge[style = {->,> = latex',pos = 0.3},color=red
			, label = $c_{1}$,labelstyle={inner sep=0pt}](6)(9);
			
			\Edge[style={bend right,->,> = latex'},color = blue, label = $c_{3}$,labelstyle={inner sep=0pt}](1)(3)
			\Edge[style={bend left,->,> = latex'},color = blue, label = $c_{3}$,labelstyle={inner sep=0pt}](8)(9)
			\Edge[style={->,> = latex'},color = red, label = $c_{1}$,labelstyle={inner sep=0pt}](4)(5)
			\end{tikzpicture}
			\caption{Step $1$.}
			\label{g:CZFS2}
		\end{subfigure}
		\medskip
		
		\begin{subfigure}{0.4\textwidth}
			\centering
			\begin{tikzpicture}[scale=0.4]
			\tikzset{VertexStyle1/.style = {shape = circle,
					ball color = white!100!black,
					text = black,
					inner sep = 1.5pt,
					outer sep = 0pt,
					minimum size = 8 pt},
				edge/.style={->,> = latex', text = black}
			}
			\tikzset{VertexStyle2/.style = {shape = circle,
					ball color = black!80!yellow,
					text = white,
					inner sep = 2pt,
					outer sep = 0pt,
					minimum size = 10 pt}}
			\node[VertexStyle2](1) at (-3,4) {$1$};
			\node[VertexStyle2](2) at (-3,0) {$2$};
			\node[VertexStyle2](3) at (-3,-4) {$3$};
			\node[VertexStyle2](4) at (3,4) {$4$};
			\node[VertexStyle2](5) at (3,0) {$5$};
			\node[VertexStyle2](6) at (3,-4) {$6$};
			\node[VertexStyle2](7) at (9,4) {$7$};
			\node[VertexStyle2](8) at (9,0) {$8$};
			\node[VertexStyle2](9) at (9,-4) {$9$};
			\Edge[style = {->,> = latex',pos = 0.5},color=green
			, label = $c_{2}$,labelstyle={inner sep=0pt}](1)(4);
			\Edge[ style = {->,> = latex',pos = 0.2},color=green
			, label = $c_{2}$,labelstyle={inner sep=0pt}](1)(5);
			\Edge[style = {->,> = latex',pos = 0.7},color=blue
			, label = $c_{3}$,labelstyle={inner sep=0pt}](1)(6);
			\Edge[ style = {->,> = latex',pos = 0.7},color=green
			, label = $c_{2}$,labelstyle={inner sep=0pt}](2)(4);
			\Edge[ style = {->,> = latex',pos = 0.3},color=red
			, label = $c_{1}$,labelstyle={inner sep=0pt}](2)(5);
			\Edge[style = {->,> = latex',pos = 0.3},color=green
			, label = $c_{2}$,labelstyle={inner sep=0pt}](3)(4);
			\Edge[style = {->,> = latex',pos = 0.5},color=blue
			, label = $c_{3}$,labelstyle={inner sep=0pt}](3)(6);
			
			\Edge[style = {->,> = latex',pos = 0.5},color=green
			, label = $c_{2}$,labelstyle={inner sep=0pt}](4)(7);
			\Edge[ style = {->,> = latex',pos = 0.2},color=green
			, label = $c_{2}$,labelstyle={inner sep=0pt}](4)(8);
			\Edge[style = {->,> = latex',pos = 0.3},color=blue
			, label = $c_{3}$,labelstyle={inner sep=0pt}](5)(8);
			\Edge[ style = {->,> = latex',pos = 0.7},color=blue
			, label = $c_{3}$,labelstyle={inner sep=0pt}](5)(9);
			\Edge[ style = {->,> = latex',pos = 0.2},color=red
			, label = $c_{1}$,labelstyle={inner sep=0pt}](6)(7);
			\Edge[style = {->,> = latex',pos = 0.3},color=red
			, label = $c_{1}$,labelstyle={inner sep=0pt}](6)(9);
			
			\Edge[style={bend right,->,> = latex'},color = blue, label = $c_{3}$,labelstyle={inner sep=0pt}](1)(3)
			\Edge[style={bend left,->,> = latex'},color = blue, label = $c_{3}$,labelstyle={inner sep=0pt}](8)(9)
			\Edge[style={->,> = latex'},color = red, label = $c_{1}$,labelstyle={inner sep=0pt}](4)(5)
			\end{tikzpicture}
			\caption{Step $2$.}
			\label{g:CZFS3}
		\end{subfigure}
		\caption{An example of a zero forcing set.}
		\label{g:CZFS}
	\end{figure}
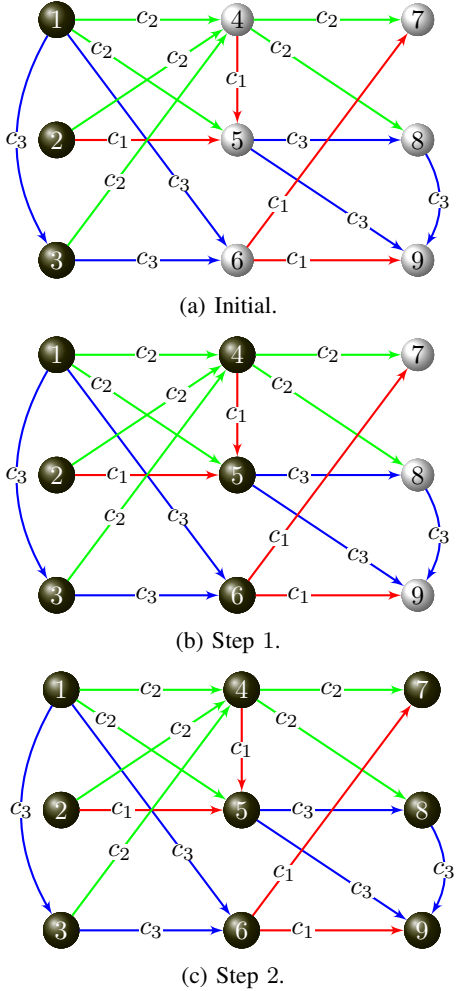
Figure \ref{g:CZFS} illustrates the repeated application of zero forcing in the context of colored graphs. 
In Figure \ref{g:CZFS1}, initially, vertices $\{1,2,3\}$ are black and the remaining vertices are white. As shown in Example~\ref{ex:coloredbg}, $\{4,5,6\}$ is a color-perfect white neighbor of $\{1,2,3\}$. Therefore, we have $\{1,2,3\} \xrightarrow{c} \{4,5,6\} $. Next, observe that the colored bipartite graph associated with $X= \{4,5,6\}$ and $Y= \{7,8,9\}$ has two perfect matchings, with identical spectrum and the same sign $1$. Hence the single equivalence class has signature $2$. As such, $ \{7,8,9\}$ is a color-perfect white neighbor of $\{4,5,6\}$. Therefore, we have $\{4,5,6\} \xrightarrow{c} \{7,8,9\}$. Consequently, we conclude that the vertex set $\{1,2,3\}$ is a zero forcing set for $\mathcal{G}(\pi)$. 
\end{ex}

%
%
Next, we explore the relationship between zero forcing sets and controllability of $(\mathcal{G}(\pi);V_{L})$.
First we show that color changes do not affect the property of controllability. This is stated in the following theorem.
\begin{thm}\label{t:cec}
	Let $\mathcal{G}(\pi)$ be a colored directed graph and let $C \subseteq V$ be a coloring set.  Suppose that $X \xrightarrow{c} Y$ with $X \subseteq C$ and $Y \subseteq V \setminus C$. Then, $(\mathcal{G} (\pi);C)$ is controllable if and only if $(\mathcal{G}(\pi);C\cup Y)$ is controllable.
\end{thm}
\begin{IEEEproof}
Due to Lemma~\ref{l:CSSCBS},  it suffices to show that $\mathcal{D}_{z}(C) = V$ if and only if 
$\mathcal{D}_{z}(C \cup Y) = V$  for all weighted graphs $\mathcal{G}(W) = (V,E,W)$ with $W \in \mathcal{W}_{\pi}(\mathcal{G})$. Here, $C$ and $C \cup Y$ are taken as zero vertex sets.

Let $W \in \mathcal{W}_{\pi}(\mathcal{G})$ and $\mathcal{G}(W) = (V,E,W)$. 
By definition of the color change rule,  $X \xrightarrow{c} Y$ means that $Y = N_{V \setminus C}(X)$ and there exists exactly one equivalence class of perfect matchings with nonzero signature in the colored bipartite graph $G = (X, Y, E_{XY}, \pi_{XY})$. 
By applying Theorem \ref{t:NoP} we then find that all matrices in the pattern class of $G$ are nonsingular.
Now, let $x_{1},x_2,\ldots,x_{n}$ be variables assigned to the vertices in $V$, with $x_{j} = 0$ for $j \in C$ and $x_j$ undetermined for the remaining vertices. 
For the vertices $j \in C$, consider the balance equations \eqref{e:BE}. 
By the fact that $W_{kj} = 0$ for all $k \in V \setminus C$ with $k \notin  N_{V \setminus C}(\{j\})$, the system of balance equations \eqref{e:BE} for the vertices  $j \in X$ can be written as
\begin{equation}\label{e:balancing}
x_{Y}^T W_{Y,X} = 0.
\end{equation}
We now observe that the submatrix $W_{Y,X}$ of $W$ belongs to the pattern class of $G$. Using the fact that all matrices in this pattern class are nonsingular, we obtain that  $x^T_{Y} = 0$. 
By the definition of the zero extension rule, we have that
$X \xrightarrow{z} Y$ for $\mathcal{G}(W)$ with the set of zero vertices $C$.
It then follows immediately that $C \cup Y \subseteq \mathcal{D}_{z}(C)$ and thus $\mathcal{D}_{z}(C \cup Y) = \mathcal{D}_{z}(C)$.
As a consequence, $C$ is a balancing set for $G(W)$ if and only if $C \cup Y$ is a balancing set for $G(W)$. Since this holds for arbitary choice of $W$ in $W_\pi(G)$, the result follows immediately from Lemma 6.

\end{IEEEproof}

By Theorem \ref{t:cec}  colored strong structural controllability is invariant under application of the color change rule. We then obtain the following corollary.
\begin{col}\label{c:cec}
	Let $\mathcal{G}(\pi)$ be a colored directed graph, let $V_L \subseteq V$ be a leader set and let $\mathcal{D}_{c}(V_{L})$ be a derived set.
	 Then $(\mathcal{G} (\pi);V_L)$ is controllable if and only if $(\mathcal{G}(\pi);\mathcal{D}_{c}(V_{L}))$ is controllable.
\end{col}

As an immediate consequence of Corollary \ref{c:cec} we arrive at the main result of this section which provides sufficient graph-theoretic condition for controllability of $(\mathcal{G}(\pi);V_{L})$.
\begin{thm}\label{t:coloredzfs}
Let $\mathcal{G}(\pi)=(V,E,\pi)$ be a colored directed graph with leader set $V_{L} \subseteq V$.
If $V_{L}$ is a zero forcing set, then $(\mathcal{G}(\pi);V_{L})$ is controllable.
\end{thm}
\begin{IEEEproof}
	The proof follows immediately from Corollary \ref{c:cec} and the fact that, trivially, $(\mathcal{G}(\pi);V)$ is controllable.
\end{IEEEproof}

To conclude this section, we will provide a counter example to show that the condition in Theorem \ref{t:coloredzfs} is not a necessary condition.
\begin{ex} \label{e:counterexample1}
	\begin{figure}[h!]
		\centering
		\begin{tikzpicture}[scale=0.5]
		\tikzset{VertexStyle1/.style = {shape = circle,
				ball color = white!100!black,
				text = black,
				inner sep = 2pt,
				outer sep = 0pt,
				minimum size = 10 pt},
			edge/.style={->,> = latex', text = black}
		}
		\tikzset{VertexStyle2/.style = {shape = circle,
				ball color = black!80!yellow,
				text = white,
				inner sep = 2pt,
				outer sep = 0pt,
				minimum size = 10 pt}}
		\node[VertexStyle2](1) at (-2,2) {$1$};
		\node[VertexStyle2](2) at (2,2) {$2$};
		\node[VertexStyle1](3) at (-5,-2) {$3$};
		\node[VertexStyle1](4) at (0,-2) {$4$};
		\node[VertexStyle1](5) at (5,-2) {$5$};
		\Edge[style={bend left,->,> = latex'},label = $c_{2}$,color = red,labelstyle={inner sep=0pt}](1)(2);
		\Edge[style={bend right,->,> = latex'},label = $c_{1}$,color = blue,labelstyle={inner sep=0pt}](1)(3);
		\Edge[ label =$c_{1}$, style = {->,> = latex',pos = 0.7},color = blue,labelstyle={inner sep=0pt}](1)(4);
		\Edge[style = {->,> = latex',pos = 0.7},color=red,label = $c_{2}$,labelstyle={inner sep=0pt}](2)(3);
		\Edge[label =$c_{2}$, style = {->,> = latex'},color=red ,labelstyle={inner sep=0pt}](2)(4);
		\Edge[style={bend left,->,> = latex'},label = $c_{1}$,color = blue,labelstyle={inner sep=0pt}](2)(5)
		\Edge[style={bend left,->,> = latex'},label = $c_{1}$,color = blue,labelstyle={inner sep=0pt}](5)(4)
		\Edge[style={bend left,->,> = latex'},label = $c_{1}$,color = blue,labelstyle={inner sep=0pt}](4)(3)
		\end{tikzpicture}
		\caption{An example to show that $V_{L}$ being a zero forcing set is not a necessary condition for controllability of $(\mathcal{G}(\pi);V_{L})$.}
		\label{g:CSDCG}
	\end{figure}
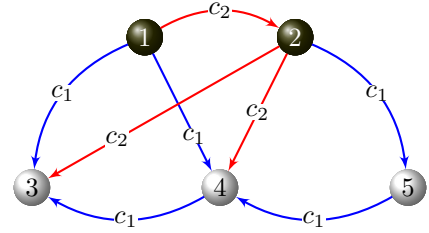
	Consider the colored  graph $\mathcal{G}(\pi)$ depicted in Figure \ref{g:CSDCG} with leader set $V_{L} = \{1,2\}$. 
	Clearly, since none of the subsets $\{1,2\}$, $\{1\}$ and $\{2\}$ have  color-perfect white neighbors, there does not exist a derived set $\mathcal{D}_{c}(V_{L})$ that equals $V$.
	Hence  $V_{L}$ is not a zero forcing set. We will show that, however, $(\mathcal{G}(\pi);V_{L})$ is controllable. 
	Due to Theorem \ref{l:CSSCBS}, it is sufficient to show that  $V_{L}$ is a balancing set for all weighted graphs $\mathcal{G}(W)$ with $W \in \mathcal{W}_{\pi}(\mathcal{G})$. 
	To do this, let $W \in \mathcal{W}_{\pi}(\mathcal{G})$ correspond to a realization $\{c_{1},c_{2}\}$ of the color set, with $c_{1}$ and $c_{2}$ nonzero real numbers.   
	Assign variables $x_{1},\ldots,x_{5}$ to the vertices in $V$. Let $x_{1} = x_{2}= 0$ and let $x_{3}$, $x_{4}$ and $x_5$ be undetermined.
	The system of balance equations \eqref{e:BE} for the vertices $1$ and $2$ in $V_L$ is then
	given by
	\begin{align}\label{e:hsbe}
	\begin{split}
	c_{1}x_{3} + c_{1}x_{4}&=0,\\
	c_{2}x_{3} + c_{2}x_{4} +c_{1}x_{5}&=0.
	\end{split}
	\end{align}
 Since $c_{1} \neq 0$ and $c_{2} \neq 0$, the homogeneous system \eqref{e:hsbe} is equivalent to the system
	\begin{align}\label{e:hsbe1}
	\begin{split}
	c_{1}x_{3} + c_{1}x_{4}&=0,\\
	c_{1}x_{5}&=0,
	\end{split}
	\end{align}
which yields $x_{5} = 0$. 
By the definition of the zero extension rule, we therefore have $\{1,2\} \xrightarrow{z} \{5\}$. 
Repeated application of the zero extension rule yields that $V_{L}$ is a balancing set. 
Since the matrix $W \in \mathcal{W}_{\pi}(\mathcal{G})$ was taken arbitrary, we conclude that $V_{L}$ is a balancing set for all weighted graphs $\mathcal{G}(W)$ with $W \in \mathcal{W}_{\pi}(\mathcal{G})$.  
Thus we have found a counter example for the necessity of the condition in Theorem \ref{t:coloredzfs}.  
\end{ex}

\section{Elementary Edge Operations and Derived Colored Graphs}

In the previous section, in Theorem \ref{t:coloredzfs}, we have established a sufficient condition
for colored strong structural controllability. In the present section we will establish another sufficient graph-theoretic condition. This new condition is based on the so-called {\em elementary edge operations}. These are operations that can be performed on the given colored graph, and that preserve colored strong structural controllability.  These edge operations on the graph are motivated by the observation that elementary operations on the systems of balance equations appearing in the zero extension rule do not modify the set of solutions to these linear equations.
Indeed, in Example~\ref{e:counterexample1}, we verified that  $\{1,2\} \xrightarrow{z} \{5\}$ for all weighted graphs $\mathcal{G}(W)$ with $W \in \mathcal{W}_{\pi}(\mathcal{G})$. 
This is due to the fact that the system of balance equations \eqref{e:hsbe} is equivalent to \eqref{e:hsbe1}, implying that $x_{5} = 0$ for all nonzero values $c_{1}$ and  $c_{2}$. 
To generalize and visualize this idea on the level of the colored graph, we now introduce the following two types of elementary edge operations.

Let $C \subseteq V$ be a coloring set, i.e., the set of vertices initially colored black. The complement $V \setminus C$ is the set of white vertices.
For two vertices $u, v \in C$ (where $u$ and $v$ can be same vertex), we define  \[\mathcal{E}_{u}(v) := \{(v,j) \in E \mid j \in N_{V \setminus C}(u)\}\]
the subset of edges between $v$ and white out-neighbors of $u$.
We now introduce the following two elementary edge operations:  
	 \begin{enumerate}
	 	\item{ ({\em Turn color}) 
	 		If all edges in $\mathcal{E}_{u}(u)$ have the same color, say $c_{i}$, then change the color of these edges to any other color in the color set.}
	 	\item{({\em Remove edges})  
	 		Assume $ N_{V \setminus C}(u) \subseteq N_{V \setminus C}(v)$. 
	 		If for any $k \in N_{V \setminus C}(u)$, the two edges $(u,k)$ and $(v,k)$ have the same color, then remove all edges in $\mathcal{E}_{u}(v)$.}
	 \end{enumerate}
The above elementary edge operations can be applied sequentially and, obviously, will not introduce new colors or add new edges. In the sequel, we will denote an edge operation by the symbol $o$.
Applying  the edge operation $o$  to $\mathcal{G}(\pi)$, we obtain a new colored  graph $\mathcal{G}'(\pi') = (V,E',\pi')$. 
We then call  $\mathcal{G}'(\pi')$ a derived graph of $\mathcal{G}(\pi)$ associated with $C$ and $o$.
We denote such derived graph by $\mathcal{G}(\pi,C,o)$. 
An application of a sequence of elementary edge operations is illustrated in the following example. 
\begin{ex}\label{ex:eeo}
	For the colored graph $\mathcal{G}(\pi) = (V,E,\pi)$ depicted in \ref{g:exeop1}, let $C = \{1,2\}$ be the coloring set. 
	For the vertex $1 \in C$,  we have $\mathcal{E}_{1}(1) = \{(1,3),(1,4)\}$ in which both edges have the same color $c_{1}$. 
	We apply the turn color operation to change the colors of $(1,3)$ and $(1,4)$ to $c_{2}$. 
	Denote this  operation by $o_{1}$.  
	We then obtain the derived colored graph $\mathcal{G}(\pi,C,o_1)$ of  $\mathcal{G}(\pi)$ with respect to $C$ and $o_1$, which is denoted by $\mathcal{G}_{1}(\pi_{1})$ and shown in \ref{g:exeop2}. 
	In addition, for the nodes $1 \mbox{ and } 2$ in $\mathcal{G}_{1}(\pi_{1})$, we have $ N_{V \setminus C}(1) \subseteq N_{V \setminus C}(2)$, where $N_{V \setminus C}(1) = \{3, 4\}$ and $N_{V \setminus C}(2) = \{3,4,5\}$. 
	Besides, for any $k \in N_{V \setminus C}(1)$, the two edges $(1,k)$ and $(2,k)$ have the same color. 
	Performing the edge removal operation denoted by $o_{2}$, we then remove all the edges in $\mathcal{E}_{1}(2) = \{(2,3),(2,4)\}$. 
	Thus we obtain the derived colored graph $\mathcal{G}_1(\pi_1, C, o_2)$ of $\mathcal{G}_1(\pi_1)$ with respect to $C$ and $o_2$, which is denoted by $\mathcal{G}_{2}(\pi_{2})$ and depicted in \ref{g:exeop3}.
\end{ex}
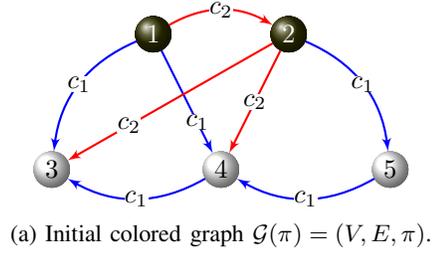
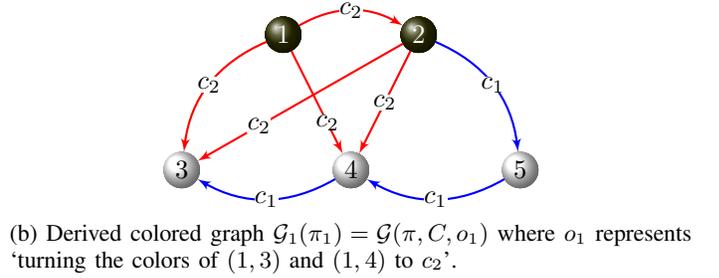
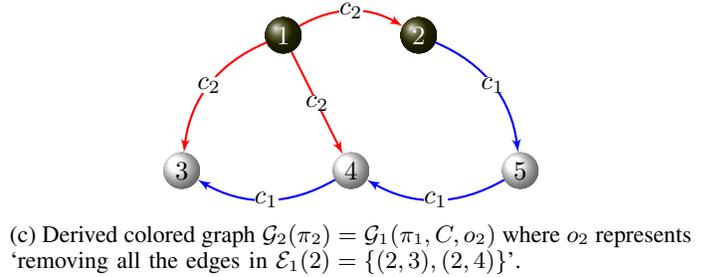
\begin{figure}[h!]
	\centering
	\begin{subfigure}{0.5\textwidth}
		\centering
		\begin{tikzpicture}[scale=0.45]
		\tikzset{VertexStyle1/.style = {shape = circle,
				ball color = white!100!black,
				text = black,
				inner sep = 2pt,
				outer sep = 0pt,
				minimum size = 10 pt},
			edge/.style={->,> = latex', text = black}
		}
		\tikzset{VertexStyle2/.style = {shape = circle,
				ball color = black!80!yellow,
				text = white,
				inner sep = 2pt,
				outer sep = 0pt,
				minimum size = 10 pt}}
		\node[VertexStyle2](1) at (-2,2) {$1$};
		\node[VertexStyle2](2) at (2,2) {$2$};
		\node[VertexStyle1](3) at (-5,-2) {$3$};
		\node[VertexStyle1](4) at (0,-2) {$4$};
		\node[VertexStyle1](5) at (5,-2) {$5$};
		\Edge[style={bend left,->,> = latex'},label = $c_{2}$,color = red,labelstyle={inner sep=0pt}](1)(2);
		\Edge[style={bend right,->,> = latex'},label = $c_{1}$,color = blue,labelstyle={inner sep=0pt}](1)(3);
		\Edge[ label =$c_{1}$, style = {->,> = latex',pos = 0.7},color = blue,labelstyle={inner sep=0pt}](1)(4);
		\Edge[style = {->,> = latex',pos = 0.7},color=red,label = $c_{2}$,labelstyle={inner sep=0pt}](2)(3);
		\Edge[label =$c_{2}$, style = {->,> = latex'},color=red ,labelstyle={inner sep=0pt}](2)(4);
		\Edge[style={bend left,->,> = latex'},label = $c_{1}$,color = blue,labelstyle={inner sep=0pt}](2)(5)
		\Edge[style={bend left,->,> = latex'},label = $c_{1}$,color = blue,labelstyle={inner sep=0pt}](5)(4)
		\Edge[style={bend left,->,> = latex'},label = $c_{1}$,color = blue,labelstyle={inner sep=0pt}](4)(3)
		\end{tikzpicture}
		\caption{Initial colored graph $\mathcal{G}(\pi) = (V,E,\pi)$.}
		\label{g:exeop1}
	\end{subfigure}
	\medskip
	
	\begin{subfigure}{0.5\textwidth}
		\centering
		\begin{tikzpicture}[scale=0.45]
		\tikzset{VertexStyle1/.style = {shape = circle,
				ball color = white!100!black,
				text = black,
				inner sep = 2pt,
				outer sep = 0pt,
				minimum size = 10 pt},
			edge/.style={->,> = latex', text = black}
		}
		\tikzset{VertexStyle2/.style = {shape = circle,
				ball color = black!80!yellow,
				text = white,
				inner sep = 2pt,
				outer sep = 0pt,
				minimum size = 10 pt}}
		\node[VertexStyle2](1) at (-2,2) {$1$};
		\node[VertexStyle2](2) at (2,2) {$2$};
		\node[VertexStyle1](3) at (-5,-2) {$3$};
		\node[VertexStyle1](4) at (0,-2) {$4$};
		\node[VertexStyle1](5) at (5,-2) {$5$};
		\Edge[style={bend left,->,> = latex'},label = $c_{2}$,color = red,labelstyle={inner sep=0pt}](1)(2);
		\Edge[style={bend right,->,> = latex'},label = $c_{2}$,color = red,labelstyle={inner sep=0pt}](1)(3);
		\Edge[ label =$c_{2}$, style = {->,> = latex',pos = 0.7},color = red,labelstyle={inner sep=0pt}](1)(4);
		\Edge[style = {->,> = latex',pos = 0.7},color = red,label = $c_{2}$,labelstyle={inner sep=0pt}](2)(3);
		\Edge[label =$c_{2}$, style = {->,> = latex'},color = red ,labelstyle={inner sep=0pt}](2)(4);
		\Edge[style={bend left,->,> = latex'},label = $c_{1}$,color = blue,labelstyle={inner sep=0pt}](2)(5)
		\Edge[style={bend left,->,> = latex'},label = $c_{1}$,color = blue,labelstyle={inner sep=0pt}](5)(4)
		\Edge[style={bend left,->,> = latex'},label = $c_{1}$,color = blue,labelstyle={inner sep=0pt}](4)(3)
		\end{tikzpicture}
		\caption{Derived colored graph $\mathcal{G}_{1}(\pi_{1}) = \mathcal{G}(\pi,C, o_1)$ where $o_1$ represents `turning the colors of $(1,3)$ and $(1,4)$ to $c_{2}$'.}
		\label{g:exeop2}
	\end{subfigure}
	\medskip
	
	\begin{subfigure}{0.5\textwidth}
		\centering
		\begin{tikzpicture}[scale=0.45]
		\tikzset{VertexStyle1/.style = {shape = circle,
				ball color = white!100!black,
				text = black,
				inner sep = 2pt,
				outer sep = 0pt,
				minimum size = 10 pt},
			edge/.style={->,> = latex', text = black}}
		\tikzset{VertexStyle2/.style = {shape = circle,
				ball color = black!80!yellow,
				text = white,
				inner sep = 2pt,
				outer sep = 0pt,
				minimum size = 10 pt}}
		\node[VertexStyle2](1) at (-2,2) {$1$};
		\node[VertexStyle2](2) at (2,2) {$2$};
		\node[VertexStyle1](3) at (-5,-2) {$3$};
		\node[VertexStyle1](4) at (0,-2) {$4$};
		\node[VertexStyle1](5) at (5,-2) {$5$};
		\Edge[style={bend right,->,> = latex'},label = $c_{2}$,color = red,labelstyle={inner sep=0pt}](1)(3);
		\Edge[style={bend left,->,> = latex'},label = $c_{2}$,color = red,labelstyle={inner sep=0pt}](1)(2);
		\Edge[ label =$c_{2}$, style = {->,> = latex',pos = 0.5},color = red,labelstyle={inner sep=0pt}](1)(4);
		\Edge[style={bend left,->,> = latex'},label = $c_{1}$,color = blue,labelstyle={inner sep=0pt}](2)(5)
		\Edge[style={bend left,->,> = latex'},label = $c_{1}$,color = blue,labelstyle={inner sep=0pt}](5)(4)
		\Edge[style={bend left,->,> = latex'},label = $c_{1}$,color = blue,labelstyle={inner sep=0pt}](4)(3)
		\end{tikzpicture}
		\caption{Derived colored graph $\mathcal{G}_{2}(\pi_{2}) = \mathcal{G}_1(\pi_1,C,o_2)$  where $o_2$ represents `removing all the edges in $\mathcal{E}_{1}(2) = \{(2,3),(2,4)\}$'.}
		\label{g:exeop3}
	\end{subfigure}	
	\caption{Example of performing elementary edge operations.
	}
	\label{g:exeop}
\end{figure}
Each elementary edge operation $o$ corresponds to a single vertex $u \in C$ or a pair of vertices $u,v \in C$. In the sequel we will denote this subset of $C$ corresponding to $o$ by $C(o)$. Thus, $C (o)$ is either a singleton or a set consisting of two elements.

Next, we study the relationship between elementary edge operations and controllability of $(\mathcal{G}(\pi);V_{L})$. First we show that elementary edge operations preserve zero extension. This issue is addressed in the following lemma.

\begin{lem}\label{l:eoet2}
	Let $\mathcal{G}(\pi)$ be a colored directed graph and $C$ be a coloring set. Let $o$ represent an edge operation and let
	$\mathcal{G'}(\pi') = \mathcal{G}(\pi,C,o)$ be a derived graph with respect to $C$ and $o$.  Let $W \in \mathcal{W}_{\pi}(\mathcal{G})$ be a weighted adjacency matrix and let $W' \in \mathcal{W}_{\pi'}(\mathcal{G}')$ be the corresponding matrix associated with the same realization of the colors. 
		Let $X \subseteq C \setminus C(o)$ and define $X' := C(o) \cup X$.
	Then, interpreting $C$ as the set of zero vertices,  for any $Y \subseteq V$ we have  $X' \xrightarrow{z} Y$  in the weighted graph $\mathcal{G}(W)$ if and only if  $X' \xrightarrow{z} Y$ in the weighted graph $\mathcal{G'}(W')$.
\end{lem}
\begin{IEEEproof}
By suitably relabeling the vertices, we may assume that $W$ has the form
\[	W = \begin{bmatrix}
	     W_{1,1} & W_{1,2} & \ldots & W_{1,6}\\
	     W_{2,1} & W_{2,2} & \ldots & W_{2,6}\\
	     W_{3,1} & W_{3,2} & \ldots & W_{3,6}\\
	     W_{4,1} & W_{4,2} & \ldots & W_{4,6}\\
	     W_{5,1} & W_{5,2} & \ldots & W_{5,6}\\
	     W_{6,1} & W_{6,2} & \ldots & W_{6,6}\\
         \end{bmatrix},
\]
where the  first row block corresponds to the vertices indexed by $ C(o)$,  the second row block corresponds to  the vertices indexed by $X$, the third row block corresponds to the vertices indexed by $C \setminus X'$, the fourth row block corresponds to the vertices indexed by $N_{V \setminus C}(C(o))$, the fifth row block corresponds to the vertices indexed by $N_{V \setminus C}(X') \setminus N_{V \setminus C}(C(o))$ and the last row block corresponds to the remaining white vertices. The column blocks of $W$ result from the same labeling. Correspondingly, the matrix $W'$ must then be equal to
\[  W' = \begin{bmatrix}
W_{1,1} & W_{1,2} & \ldots & W_{1,6}\\
W_{2,1} & W_{2,2} & \ldots & W_{2,6}\\
W_{3,1} & W_{3,2} & \ldots & W_{3,6}\\
W'_{4,1} & W_{4,2} & \ldots & W_{4,6}\\
W_{5,1} & W_{5,2} & \ldots & W_{5,6}\\
W_{6,1} & W_{6,2} & \ldots & W_{6,6}\\
\end{bmatrix}.
\]
for some matrix $W'_{4,1}$.
Since the fourth and fifth row blocks correspond to the vertices indexed by $N_{V \setminus C}(C(o)) $ and $N_{V \setminus C}(X') \setminus N_{V \setminus C}(C(o))$, respectively, it follows easily that $W_{5,1} = \mathbf{0}$, $W_{6,1} = \mathbf{0}$ and $W_{6,2} = \mathbf{0}$.
Consider the submatrices $W_{N_{V \setminus C}(X'), X'} = \begin{bmatrix}
W_{4,1} & W_{4,2} \\
\mathbf{0} & W_{5,2} \\
\end{bmatrix}$ and $W'_{N_{V \setminus C}(X'), X'}= \begin{bmatrix}
W'_{4,1} & W_{4,2} \\
\mathbf{0} & W_{5,2} \\
\end{bmatrix}$  of $W$ and $W'$, respectively.  
We then distinguish two cases:
\begin{enumerate}
	\item Suppose the edge operation $o$ represents a color turn operation. In that case, $C(o)$ only contains one vertex, in other words, both $W_{4,1}$ and $W'_{4,1}$ consist of only one column. Hence, it follows that $W'_{4,1} = \alpha W_{4,1}$ for a suitable nonzero real number $\alpha$.
	\item Suppose the edge operation $o$ represents an edge removal operation. In that case $C(o)$ contains two vertices, say $u$ and $v$, and both $W_{4,1}$ and $W'_{4,1}$ consist of two columns. We may assume that $u$ and $v$ correspond to the first and second column of these matrices, respectively, and the edges in $\mathcal{E}_u(v)$ are removed. This implies that $$W'_{4,1} = W_{4,1} \begin{bmatrix}
	1 & -1 \\
	0 & 1 \\
	\end{bmatrix}.$$
\end{enumerate}
Clearly, $W_{N_{V \setminus C}(X'), X'}$ and $W'_{N_{V \setminus C}(X'), X'}$ are column equivalent. 
Next, again assign variables $x_{1},\ldots,x_{n}$ to every vertex in $V$, where  $x_{i}$ is equal to $0$ if $i \in C$ and otherwise undetermined.  
For the vertex $j \in C$ we consider the balance equation \eqref{e:BE}.
By the fact that $W_{kj} = 0$ for all $k \in V \setminus C$ with $k \notin  N_{V \setminus C}(\{j\})$ and $N_{V \setminus C}(\{j\}) \subseteq N_{V \setminus C}(X')$, equation \eqref{e:BE} is equivalent to
\begin{equation}\label{e:BE2}
	\sum_{k \in N_{V \setminus C}(X')} x_{k}W_{kj} = 0.
\end{equation}
Again using the notation for the submatrix $W_{N_{V \setminus C}(X'), X'}$ and subvector $x_{N_{V \setminus C}(X')}$, we can rewrite the system of balance equations \eqref{e:BE2} for $j \in X'$ as
\begin{equation}\label{e:ghe1}
    x_{N_{V \setminus C}(X')}^{T}W_{N_{V \setminus C}(X'), X'} = 0.
\end{equation}
Similarly, for the graph $\mathcal{G}' (W')$, we obtain the following system of balance equations for $j \in X'$:
\begin{equation}\label{e:ghe2}
      x_{N_{V \setminus C}(X')}^{T}W'_{N_{V \setminus C}(X'), X'} = 0.
\end{equation}
Since $W'_{N_{V \setminus C}(X'), X'}$ and $W_{N_{V \setminus C}(X'), X'}$ are column equivalent, the solution sets of \eqref{e:ghe1} and \eqref{e:ghe2} coincide.
By definition of the zero extension rule we therefore have that, for any vertex set $Y$,  $X' \xrightarrow{z} Y$ in $\mathcal{G}(W)$ if and only if $X' \xrightarrow{z} Y$ in $\mathcal{G}' (W')$. 
This completes the proof.		
\end{IEEEproof}

It follows from the previous that colored strong structural controllability is preserved under elementary edge operations. Indeed, we have
\begin{thm}\label{t:eeo}
	Let $\mathcal{G}(\pi)$ be a colored directed graph, $V_L \subseteq V$ be a leader set, and $o$ an elementary edge operation. Let $\mathcal{G}' (\pi') = \mathcal{G}(\pi,V_L,o)$ be a derived colored graph of $\mathcal{G}(\pi)$ with respect to $V_L$ and $o$. 
	Then we have that $(\mathcal{G}(\pi);V_{L})$ is controllable if and only if $(\mathcal{G'} (\pi');V_L)$ is controllable.
\end{thm}
\begin{IEEEproof}
	The proof follows from Lemma \ref{l:CSSCBS} and Lemma \ref{l:eoet2}.
\end{IEEEproof}

As an immediate consequence of Theorem \ref{t:eeo} and Theorem \ref{t:coloredzfs} we see that
if the leader set $V_L$ of the original colored graph $\mathcal{G}(\pi)$ is a zero forcing set for the derived graph $\mathcal{G}' (\pi') = \mathcal{G}(\pi,V_L,o)$, then $(\mathcal{G} (\pi);V_L)$ is controllable. 
Obviously, this result can immediately be extended to derived graphs obtained by applying a finite sequence of edge operations. This leads to the following sufficient graph theoretic condition for controllability of 
$(\mathcal{G}(\pi);V_{L})$.
\begin{col} \label{c:ecs}
Let $\mathcal{G}(\pi)$ be a colored directed graph and let $V_L$ be a leader set. Let $\mathcal{G}'(\pi')$ be a colored graph obtained by applying finitely many elementary edge operations. Then $(\mathcal{G}(\pi),V_L)$ is controllable if $V_L$ is a zero forcing set for $\mathcal{G}'(\pi')$.
\end{col}
\begin{ex}
	Again consider the colored graph in Example \ref{e:counterexample1}. We already saw that $V_L =\{1,2\}$ is not a zero forcing set. However, we also showed that we do have strong structural controllability for this colored graph. This can now also be shown graph theoretically by means of Corollary \ref{c:ecs}: the leader set $V_L$ is a zero forcing set for the derived graph in Figure \ref{g:exeop3}, so the original colored graph in Figure \ref{g:exeop1} yields a controllable system.
\end{ex}

By combining Theorem \ref{t:eeo} and Corollary \ref{c:cec} we are now in the position to establish yet another procedure for checking controllability of a given colored graph $(\mathcal{G} (\pi);V_L)$. First, distinguish the following two steps:
	\begin{enumerate}
		\item
		As the first step, apply the color change operation to compute a derived set $\mathcal{D}_{c}(V_{L})$. If this derived set is equal to $V$ we have controllability. If not, we can not yet decide whether we have controllability or not. 
		\item
		As a next step, then, apply an edge operation $o$ to $\mathcal{G} (\pi)$ to obtain $\mathcal{G}_{1}(\pi_{1})$, where $\mathcal{G}_{1}(\pi_{1}) = \mathcal{G} (\pi,\mathcal{D}_{c}(V_{L}), o)$ is a derived graph of $\mathcal{G}(\pi)$ with coloring set $\mathcal{D}_{c}(V_{L})$ and edge operation $o$. 
	\end{enumerate}
	By Theorem \ref{t:eeo} and Corollary \ref{c:cec}, it is straightforward to verify that $(\mathcal{G} (\pi);V_L)$ is controllable if and only if $(\mathcal{G}_{1}(\pi_{1});\mathcal{D}_{c}(V_{L}))$ is controllable. 

We can now repeat steps 1 and 2, applying them to $\mathcal{G}_{1}(\pi_{1})$.
	Successive and alternating application of these two steps transforms the original leader set $V_L$ 
	using several color change operations associated with the several derived graphs appearing in the process. After finitely many iterations we thus arrive at a so called {\em edge-operations-color-change derived set} of $V_L$, that will be denoted by $\mathcal{D}_{ec}(C)$. This set will remain unchanged in case we again apply step 1 or step 2. Since controllability is preserved, we arrive at the following theorem.
\begin{thm}\label{t:eocd}
	Let $\mathcal{G}(\pi)$ be a colored directed graph and let $V_L \subseteq V$ be a leader set. Let $\mathcal{D}_{ec}(V_{L})$ be an edge-operations-color-change derived set of $V_{L}$.
	We then have that $(\mathcal{G} (\pi);V_L)$ is controllable if $\mathcal{D}_{ec}(V_{L}) = V$.
\end{thm}
\begin{remark}
	Obviously, a derived set $\mathcal{D}_{c}(V_{L})$ of $V_{L}$ in $\mathcal{G}(\pi)$ is always contained in  an edge-operations-color-change derived set $\mathcal{D}_{ec}(V_{L})$ of $V_{L}$.
		Hence the condition in Theorem \ref{t:eocd} is weaker than the conditions in Theorem \ref{t:coloredzfs} and Corollary \ref{c:ecs}.
\end{remark}

In the following example we illustrate the application of Theorem \ref{t:eocd} to check controllability of a given colored graph and leader set.
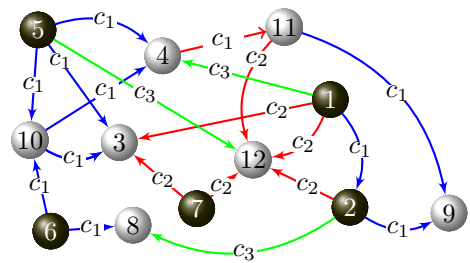
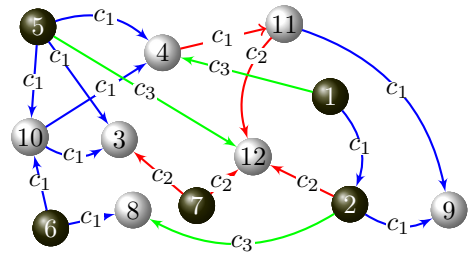
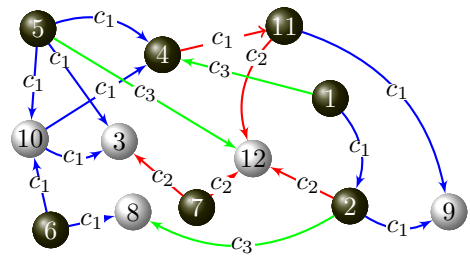
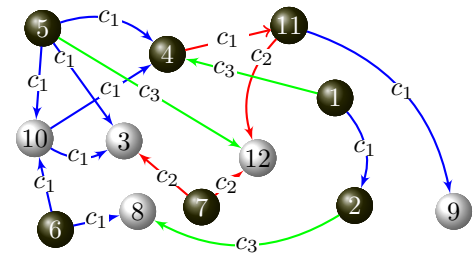
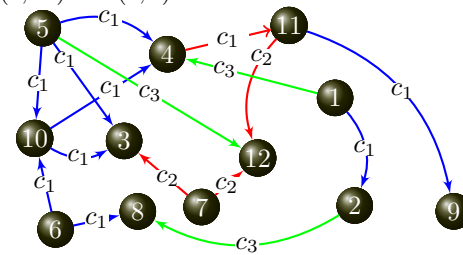
\begin{figure}[h!]\label{fi4al1}
	\centering
	\begin{subfigure}{0.45\textwidth}
		\centering
		\begin{tikzpicture}[scale=0.16]
		\tikzset{VertexStyle1/.style = {shape = circle,
				ball color = white!100!black,
				text = black,
				inner sep = 2pt,
				outer sep = 0pt,
				minimum size = 8 pt}}
		\tikzset{VertexStyle2/.style = {shape = circle,
				ball color = black!80!yellow,
				text = white,
				inner sep = 2pt,
				outer sep = 0pt,
				minimum size = 8 pt}}
		\tikzset{VertexStyle3/.style = {shape = circle,
				ball color =  white!100!black,
				text = black,
				inner sep = 0.8pt,
				outer sep = 0pt,
				minimum size = 8 pt}}
		\node[VertexStyle2](1) at (8.840,2.3005) {$1$};
		\node[VertexStyle2](2) at (10.4707,-6.6488) {$2$};
		\node[VertexStyle1](3) at (-8.6887,-1.2771) {$3$};
		\node[VertexStyle1](4) at (-5.0950,5.9327) {$4$};
		\node[VertexStyle2](5) at (-15.4428,8.0927) {$5$};
		\node[VertexStyle2](6) at (-14.3838,-8.6365) {$6$};
		\node[VertexStyle2](7) at (-2.2519 , -6.7736) {$7$};
		\node[VertexStyle1](8) at (-7.5493,-8.1412) {$8$};
		\node[VertexStyle1](9) at (18.7030,-7.1350) {$9$};
		\node[VertexStyle3](10) at (-16.1152, -1.0685) {$10$};
		\node[VertexStyle3](11) at (5.0224, 8.3263) {$11$};
		\node[VertexStyle3](12) at (2.4145,  -2.6967) {$12$};
		\Edge[ style = {->,> = latex'},color=blue,label = $c_{1}$ ,labelstyle={inner sep=0pt}](6)(8);
		\Edge[ style = {->,> = latex'},color=blue,label = $c_{1}$ ,labelstyle={inner sep=0pt}](6)(10);
		\Edge[ style = {->,> = latex'},color=blue,label = $c_{1}$ ,labelstyle={inner sep=0pt}](5)(10);
		\Edge[ style = {bend right,->,> = latex'},color=blue,label = $c_{1}$ ,labelstyle={inner sep=0pt}](10)(3);
		\Edge[ style = {->,> = latex',pos = 0.2},color=blue,label = $c_{1}$ ,labelstyle={inner sep=0pt}](5)(3);
		\Edge[ style = {bend right,->,> = latex'},color=blue,label = $c_{1}$ ,labelstyle={inner sep=0pt}](10)(3);
		\Edge[ style = ->,color=red ,label = $c_{1}$](4)(11);
		\Edge[style={bend left,->,> = latex'},label = $c_{2}$,color = red,labelstyle={inner sep=0pt}](1)(12);
		\Edge[style={bend left,->,> = latex'},label = $c_{1}$,color = blue,labelstyle={inner sep=0pt}](1)(2);
		\Edge[ label =$c_{2}$, style = {->,> = latex',pos = 0.2},color = red,labelstyle={inner sep=0pt}](1)(3);
		\Edge[style = {->,> = latex',pos = 0.4},color=red,label = $c_{2}$,labelstyle={inner sep=0pt}](7)(3);
		\Edge[style = {->,> = latex',pos = 0.4},color=red,label = $c_{2}$,labelstyle={inner sep=0pt}](7)(12);
		\Edge[label =$c_{2}$, style = {->,> = latex',pos = 0.4},color=red ,labelstyle={inner sep=0pt}](2)(12);
		\Edge[style={bend right,->,> = latex'},label = $c_{1}$,color = blue,labelstyle={inner sep=0pt}](2)(9)
		\Edge[style={bend left,->,> = latex'},label = $c_{1}$,color = blue,labelstyle={inner sep=0pt}](5)(4)
		\Edge[style={bend left,->,> = latex'},label = $c_{1}$,color = blue,labelstyle={inner sep=0pt}](11)(9)
		\Edge[style={bend right,->,> = latex',pos = 0.2},label = $c_{2}$,color = red,labelstyle={inner sep=0pt}](11)(12)
		\Edge[style={->,> = latex',pos = 0.6},label = $c_{1}$,color = blue,labelstyle={inner sep=0pt}](10)(4)
		\Edge[ style = {->,> = latex',pos = 0.5},color=green,label = $c_{3}$ ,labelstyle={inner sep=0pt}](5)(12)
		\Edge[ label =$c_{3}$, style = {->,> = latex',pos = 0.7},color = green,labelstyle={inner sep=0pt}](1)(4)
		\Edge[style={bend left,->,> = latex'},label = $c_{3}$,color = green,labelstyle={inner sep=0pt}](2)(8)
		\end{tikzpicture}
		\caption{Initial colored graph $\mathcal{G}(\pi) = (V,E,\pi)$ with coloring set $V_L  = \{1,2,5,6,7\}$. Let $\mathcal{G}_{0}(\pi_{0}) = \mathcal{G}(\pi).$ Compute a derived set $\mathcal{D}_{c}(V_L) = \{1,2,5,6,7\}$ of $V_L$ in $\mathcal{G}_{0}(\pi_{0})$ and set $\mathcal{D}_{0} = \mathcal{D}_{c}(V_L)$.}
		\label{g:exa1}
	\end{subfigure}
	\smallskip
	
	\begin{subfigure}{0.45\textwidth}
		\centering
		\begin{tikzpicture}[scale=0.16]
		\tikzset{VertexStyle1/.style = {shape = circle,
				ball color = white!100!black,
				text = black,
				inner sep = 2pt,
				outer sep = 0pt,
				minimum size = 8 pt}}
		\tikzset{VertexStyle2/.style = {shape = circle,
				ball color = black!80!yellow,
				text = white,
				inner sep = 2pt,
				outer sep = 0pt,
				minimum size = 8 pt}}
		\tikzset{VertexStyle3/.style = {shape = circle,
				ball color =  white!100!black,
				text = black,
				inner sep = 0.8pt,
				outer sep = 0pt,
				minimum size = 8 pt}}
		\tikzset{VertexStyle4/.style = {shape = circle,
				ball color = black!80!yellow,
				text = white,
				inner sep = 0.8pt,
				outer sep = 0pt,
				minimum size = 8 pt}}
		\node[VertexStyle2](1) at (8.840,2.3005) {$1$};
		\node[VertexStyle2](2) at (10.4707,-6.6488) {$2$};
		\node[VertexStyle1](3) at (-8.6887,-1.2771) {$3$};
		\node[VertexStyle1](4) at (-5.0950,5.9327) {$4$};
		\node[VertexStyle2](5) at (-15.4428,8.0927) {$5$};
		\node[VertexStyle2](6) at (-14.3838,-8.6365) {$6$};
		\node[VertexStyle2](7) at (-2.2519 , -6.7736) {$7$};
		\node[VertexStyle1](8) at (-7.5493,-7.1412) {$8$};
		\node[VertexStyle1](9) at (18.7030,-7.1350) {$9$};
		\node[VertexStyle3](10) at (-16.1152, -1.0685) {$10$};
		\node[VertexStyle3](11) at (5.0224, 8.3263) {$11$};
		\node[VertexStyle3](12) at (2.4145,  -2.6967) {$12$};
		\Edge[ style = {->,> = latex'},color=blue,label = $c_{1}$ ,labelstyle={inner sep=0pt}](6)(8);
		\Edge[ style = {->,> = latex'},color=blue,label = $c_{1}$ ,labelstyle={inner sep=0pt}](6)(10);
		\Edge[ style = {->,> = latex'},color=blue,label = $c_{1}$ ,labelstyle={inner sep=0pt}](5)(10);
		\Edge[ style = {bend right,->,> = latex'},color=blue,label = $c_{1}$ ,labelstyle={inner sep=0pt}](10)(3);
		\Edge[ style = {->,> = latex',pos = 0.2},color=blue,label = $c_{1}$ ,labelstyle={inner sep=0pt}](5)(3);
		\Edge[ style = ->,color=red ,label = $c_{1}$](4)(11);
		\Edge[style={bend left,->,> = latex'},label = $c_{1}$,color = blue,labelstyle={inner sep=0pt}](1)(2);
		\Edge[style = {->,> = latex',pos = 0.4},color=red,label = $c_{2}$,labelstyle={inner sep=0pt}](7)(3);
		\Edge[style = {->,> = latex',pos = 0.4},color=red,label = $c_{2}$,labelstyle={inner sep=0pt}](7)(12);
		\Edge[label =$c_{2}$, style = {->,> = latex',pos = 0.4},color=red ,labelstyle={inner sep=0pt}](2)(12);
		\Edge[style={bend right,->,> = latex'},label = $c_{1}$,color = blue,labelstyle={inner sep=0pt}](2)(9)
		\Edge[style={bend left,->,> = latex'},label = $c_{1}$,color = blue,labelstyle={inner sep=0pt}](5)(4)
		\Edge[style={bend left,->,> = latex'},label = $c_{1}$,color = blue,labelstyle={inner sep=0pt}](11)(9)
		\Edge[style={bend right,->,> = latex',pos = 0.2},label = $c_{2}$,color = red,labelstyle={inner sep=0pt}](11)(12)
		\Edge[style={->,> = latex',pos = 0.6},label = $c_{1}$,color = blue,labelstyle={inner sep=0pt}](10)(4)
		\Edge[ style = {->,> = latex',pos = 0.5},color=green,label = $c_{3}$ ,labelstyle={inner sep=0pt}](5)(12)
		\Edge[ label =$c_{3}$, style = {->,> = latex',pos = 0.7},color = green,labelstyle={inner sep=0pt}](1)(4)
		\Edge[style={bend left,->,> = latex'},label = $c_{3}$,color = green,labelstyle={inner sep=0pt}](2)(8)
		\end{tikzpicture}
		\caption{ Derived colored graph $\mathcal{G}_{1}(\pi_{1}) = \mathcal{G}(\pi,\mathcal{D}_{0},o_0)$ of $\mathcal{G}(\pi)$ with respect to $\mathcal{D}_{0}$ and  $o_0$ such that $o_0$ represents `removing edges $(1,12) \mbox{ and } (1,3)$'.}
		\label{g:exa2}
	\end{subfigure}
	\smallskip
	
	\begin{subfigure}{0.45\textwidth}
		\centering
		\begin{tikzpicture}[scale=0.16]
		\tikzset{VertexStyle1/.style = {shape = circle,
				ball color = white!100!black,
				text = black,
				inner sep = 2pt,
				outer sep = 0pt,
				minimum size = 8 pt}}
		\tikzset{VertexStyle2/.style = {shape = circle,
				ball color = black!80!yellow,
				text = white,
				inner sep = 2pt,
				outer sep = 0pt,
				minimum size = 8 pt}}
		\tikzset{VertexStyle3/.style = {shape = circle,
				ball color =  white!100!black,
				text = black,
				inner sep = 0.8pt,
				outer sep = 0pt,
				minimum size = 8 pt}}
		\tikzset{VertexStyle4/.style = {shape = circle,
				ball color = black!80!yellow,
				text = white,
				inner sep = 0.8pt,
				outer sep = 0pt,
				minimum size = 8 pt}}
		\node[VertexStyle2](1) at (8.840,2.3005) {$1$};
		\node[VertexStyle2](2) at (10.4707,-6.6488) {$2$};
		\node[VertexStyle1](3) at (-8.6887,-1.2771) {$3$};
		\node[VertexStyle2](4) at (-5.0950,5.9327) {$4$};
		\node[VertexStyle2](5) at (-15.4428,8.0927) {$5$};
		\node[VertexStyle2](6) at (-14.3838,-8.6365) {$6$};
		\node[VertexStyle2](7) at (-2.2519 , -6.7736) {$7$};
		\node[VertexStyle1](8) at (-7.5493,-7.1412) {$8$};
		\node[VertexStyle1](9) at (18.7030,-7.1350) {$9$};
		\node[VertexStyle3](10) at (-16.1152, -1.0685) {$10$};
		\node[VertexStyle4](11) at (5.0224, 8.3263) {$11$};
		\node[VertexStyle3](12) at (2.4145,  -2.6967) {$12$};
		\Edge[ style = {->,> = latex'},color=blue,label = $c_{1}$ ,labelstyle={inner sep=0pt}](6)(8);
		\Edge[ style = {->,> = latex'},color=blue,label = $c_{1}$ ,labelstyle={inner sep=0pt}](6)(10);
		\Edge[ style = {->,> = latex'},color=blue,label = $c_{1}$ ,labelstyle={inner sep=0pt}](5)(10);
		\Edge[ style = {bend right,->,> = latex'},color=blue,label = $c_{1}$ ,labelstyle={inner sep=0pt}](10)(3);
		\Edge[ style = {->,> = latex',pos = 0.2},color=blue,label = $c_{1}$ ,labelstyle={inner sep=0pt}](5)(3);
		\Edge[ style = ->,color=red ,label = $c_{1}$](4)(11);
		\Edge[style={bend left,->,> = latex'},label = $c_{1}$,color = blue,labelstyle={inner sep=0pt}](1)(2);
		
		\Edge[style = {->,> = latex',pos = 0.4},color=red,label = $c_{2}$,labelstyle={inner sep=0pt}](7)(3);
		\Edge[style = {->,> = latex',pos = 0.4},color=red,label = $c_{2}$,labelstyle={inner sep=0pt}](7)(12);
		\Edge[label =$c_{2}$, style = {->,> = latex',pos = 0.4},color=red ,labelstyle={inner sep=0pt}](2)(12);
		
		\Edge[style={bend right,->,> = latex'},label = $c_{1}$,color = blue,labelstyle={inner sep=0pt}](2)(9)
		\Edge[style={bend left,->,> = latex'},label = $c_{1}$,color = blue,labelstyle={inner sep=0pt}](5)(4)
		\Edge[style={bend left,->,> = latex'},label = $c_{1}$,color = blue,labelstyle={inner sep=0pt}](11)(9)
		\Edge[style={bend right,->,> = latex',pos = 0.2},label = $c_{2}$,color = red,labelstyle={inner sep=0pt}](11)(12)
		\Edge[style={->,> = latex',pos = 0.6},label = $c_{1}$,color = blue,labelstyle={inner sep=0pt}](10)(4)
		\Edge[ style = {->,> = latex',pos = 0.5},color=green,label = $c_{3}$ ,labelstyle={inner sep=0pt}](5)(12)
		\Edge[ label =$c_{3}$, style = {->,> = latex',pos = 0.7},color = green,labelstyle={inner sep=0pt}](1)(4)
		\Edge[style={bend left,->,> = latex'},label = $c_{3}$,color = green,labelstyle={inner sep=0pt}](2)(8)
		\end{tikzpicture}
		\caption{ Derived set $\mathcal{D}_{1} = \{1,2,4,5,6,7,11\}$ of $\mathcal{D}_{0}$ in the colored graph $\mathcal{G}_{1}(\pi_{1})$.}
		\label{g:exa3}
	\end{subfigure}
	\smallskip
	
	\begin{subfigure}{0.45\textwidth}
		\centering
		\begin{tikzpicture}[scale=0.16]
		\tikzset{VertexStyle1/.style = {shape = circle,
				ball color = white!100!black,
				text = black,
				inner sep = 2pt,
				outer sep = 0pt,
				minimum size = 8 pt}}
		\tikzset{VertexStyle2/.style = {shape = circle,
				ball color = black!80!yellow,
				text = white,
				inner sep = 2pt,
				outer sep = 0pt,
				minimum size = 8 pt}}
		\tikzset{VertexStyle3/.style = {shape = circle,
				ball color =  white!100!black,
				text = black,
				inner sep = 0.8pt,
				outer sep = 0pt,
				minimum size = 8 pt}}
		\tikzset{VertexStyle4/.style = {shape = circle,
				ball color = black!80!yellow,
				text = white,
				inner sep = 0.8pt,
				outer sep = 0pt,
				minimum size = 8 pt}}
		\node[VertexStyle2](1) at (8.840,2.3005) {$1$};
		\node[VertexStyle2](2) at (10.4707,-6.6488) {$2$};
		\node[VertexStyle1](3) at (-8.6887,-1.2771) {$3$};
		\node[VertexStyle2](4) at (-5.0950,5.9327) {$4$};
		\node[VertexStyle2](5) at (-15.4428,8.0927) {$5$};
		\node[VertexStyle2](6) at (-14.3838,-8.6365) {$6$};
		\node[VertexStyle2](7) at (-2.2519 , -6.7736) {$7$};
		\node[VertexStyle1](8) at (-7.5493,-7.1412) {$8$};
		\node[VertexStyle1](9) at (18.7030,-7.1350) {$9$};
		\node[VertexStyle3](10) at (-16.1152, -1.0685) {$10$};
		\node[VertexStyle4](11) at (5.0224, 8.3263) {$11$};
		\node[VertexStyle3](12) at (2.4145,  -2.6967) {$12$};
		\Edge[ style = {->,> = latex'},color=blue,label = $c_{1}$ ,labelstyle={inner sep=0pt}](6)(8);
		\Edge[ style = {->,> = latex'},color=blue,label = $c_{1}$ ,labelstyle={inner sep=0pt}](6)(10);
		\Edge[ style = {->,> = latex'},color=blue,label = $c_{1}$ ,labelstyle={inner sep=0pt}](5)(10);
		\Edge[ style = {bend right,->,> = latex'},color=blue,label = $c_{1}$ ,labelstyle={inner sep=0pt}](10)(3);
		\Edge[ style = {->,> = latex',pos = 0.2},color=blue,label = $c_{1}$ ,labelstyle={inner sep=0pt}](5)(3);
		\Edge[ style = ->,color=red ,label = $c_{1}$](4)(11);
		\Edge[style={bend left,->,> = latex'},label = $c_{1}$,color = blue,labelstyle={inner sep=0pt}](1)(2);
		
		\Edge[style = {->,> = latex',pos = 0.4},color=red,label = $c_{2}$,labelstyle={inner sep=0pt}](7)(3);
		\Edge[style = {->,> = latex',pos = 0.4},color=red,label = $c_{2}$,labelstyle={inner sep=0pt}](7)(12);
		
		\Edge[style={bend left,->,> = latex'},label = $c_{1}$,color = blue,labelstyle={inner sep=0pt}](5)(4)
		\Edge[style={bend left,->,> = latex'},label = $c_{1}$,color = blue,labelstyle={inner sep=0pt}](11)(9)
		\Edge[style={bend right,->,> = latex',pos = 0.2},label = $c_{2}$,color = red,labelstyle={inner sep=0pt}](11)(12)
		\Edge[style={->,> = latex',pos = 0.6},label = $c_{1}$,color = blue,labelstyle={inner sep=0pt}](10)(4)
		\Edge[ style = {->,> = latex',pos = 0.5},color=green,label = $c_{3}$ ,labelstyle={inner sep=0pt}](5)(12)
		\Edge[ label =$c_{3}$, style = {->,> = latex',pos = 0.7},color = green,labelstyle={inner sep=0pt}](1)(4)
		\Edge[style={bend left,->,> = latex'},label = $c_{3}$,color = green,labelstyle={inner sep=0pt}](2)(8)
		\end{tikzpicture}
		\caption{Derived colored graph $\mathcal{G}_{2}(\pi_{2}) = \mathcal{G}_{1}(\pi_{1},\mathcal{D}_{1},o_1)$ with $\mathcal{D}_{1} = \{1,2,4,5,6,7,11\}$ and $o_1$ such that $o_1$ represents `removing edges $(2,12) \mbox{ and } (2,9)$'.}
		\label{g:exa4}
	\end{subfigure}
	
	\begin{subfigure}{0.45\textwidth}
		\centering
		\begin{tikzpicture}[scale=0.16]
		\tikzset{VertexStyle1/.style = {shape = circle,
				ball color = white!100!black,
				text = black,
				inner sep = 2pt,
				outer sep = 0pt,
				minimum size = 8 pt}}
		\tikzset{VertexStyle2/.style = {shape = circle,
				ball color = black!80!yellow,
				text = white,
				inner sep = 2pt,
				outer sep = 0pt,
				minimum size = 8 pt}}
		\tikzset{VertexStyle3/.style = {shape = circle,
				ball color =  white!100!black,
				text = black,
				inner sep = 0.8pt,
				outer sep = 0pt,
				minimum size = 8 pt}}
		\tikzset{VertexStyle4/.style = {shape = circle,
				ball color = black!80!yellow,
				text = white,
				inner sep = 0.8pt,
				outer sep = 0pt,
				minimum size = 8 pt}}
		\node[VertexStyle2](1) at (8.840,2.3005) {$1$};
		\node[VertexStyle2](2) at (10.4707,-6.6488) {$2$};
		\node[VertexStyle2](3) at (-8.6887,-1.2771) {$3$};
		\node[VertexStyle2](4) at (-5.0950,5.9327) {$4$};
		\node[VertexStyle2](5) at (-15.4428,8.0927) {$5$};
		\node[VertexStyle2](6) at (-14.3838,-8.6365) {$6$};
		\node[VertexStyle2](7) at (-2.2519 , -6.7736) {$7$};
		\node[VertexStyle2](8) at (-7.5493,-7.1412) {$8$};
		\node[VertexStyle2](9) at (18.7030,-7.1350) {$9$};
		\node[VertexStyle4](10) at (-16.1152, -1.0685) {$10$};
		\node[VertexStyle4](11) at (5.0224, 8.3263) {$11$};
		\node[VertexStyle4](12) at (2.4145,  -2.6967) {$12$};
		\Edge[ style = {->,> = latex'},color=blue,label = $c_{1}$ ,labelstyle={inner sep=0pt}](6)(8);
		\Edge[ style = {->,> = latex'},color=blue,label = $c_{1}$ ,labelstyle={inner sep=0pt}](6)(10);
		\Edge[ style = {->,> = latex'},color=blue,label = $c_{1}$ ,labelstyle={inner sep=0pt}](5)(10);
		\Edge[ style = {bend right,->,> = latex'},color=blue,label = $c_{1}$ ,labelstyle={inner sep=0pt}](10)(3);
		\Edge[ style = {->,> = latex',pos = 0.2},color=blue,label = $c_{1}$ ,labelstyle={inner sep=0pt}](5)(3);
		\Edge[ style = ->,color=red ,label = $c_{1}$](4)(11);
		\Edge[style={bend left,->,> = latex'},label = $c_{1}$,color = blue,labelstyle={inner sep=0pt}](1)(2);
		
		\Edge[style = {->,> = latex',pos = 0.4},color=red,label = $c_{2}$,labelstyle={inner sep=0pt}](7)(3);
		\Edge[style = {->,> = latex',pos = 0.4},color=red,label = $c_{2}$,labelstyle={inner sep=0pt}](7)(12);
		
		\Edge[style={bend left,->,> = latex'},label = $c_{1}$,color = blue,labelstyle={inner sep=0pt}](5)(4)
		\Edge[style={bend left,->,> = latex'},label = $c_{1}$,color = blue,labelstyle={inner sep=0pt}](11)(9)
		\Edge[style={bend right,->,> = latex',pos = 0.2},label = $c_{2}$,color = red,labelstyle={inner sep=0pt}](11)(12)
		\Edge[style={->,> = latex',pos = 0.6},label = $c_{1}$,color = blue,labelstyle={inner sep=0pt}](10)(4)
		\Edge[ style = {->,> = latex',pos = 0.5},color=green,label = $c_{3}$ ,labelstyle={inner sep=0pt}](5)(12)
		\Edge[ label =$c_{3}$, style = {->,> = latex',pos = 0.7},color = green,labelstyle={inner sep=0pt}](1)(4)
		\Edge[style={bend left,->,> = latex'},label = $c_{3}$,color = green,labelstyle={inner sep=0pt}](2)(8)
		\end{tikzpicture}
		\caption{Derived set $\mathcal{D}_{2} = V$ of $\mathcal{D}_{1}$ in the colored graph $\mathcal{G}_{2}(\pi_{2})$. Return that $(\mathcal{G}(\pi);V_{L})$ is controllable.}
		\label{g:exa5}
	\end{subfigure}	
	\caption{An example of application of Theorem \ref{t:eocd}}
	\label{g:exa}
\end{figure}
\begin{ex}\label{ex:algorithm}
	Consider the colored graph $\mathcal{G}(\pi) = (V,E,\pi)$ depicted in Figure \ref{g:exa1} with  $V_L = \{1,2,5,6,7\}$ the leader set.  
	To start with, we compute a derived set $\mathcal{D}_{c}(V_L) = \{1,2,5,6,7\}$ of $V_L$ in $\mathcal{G}(\pi)$, and denote it  by $\mathcal{D}_{0}$. 
	For the vertices $1 , 7 \in \mathcal{D}_{0}$, in $\mathcal{G}(\pi)$ we have $ N_{V \setminus \mathcal{D}_{0}}(7) \subseteq N_{V \setminus \mathcal{D}_{0}}(1)$, and  for any $k \in N_{V \setminus \mathcal{D}_{0}}(7)$, the two edges $(1,k)$ and $(7,k)$ have the same color. Thus we remove all edges in $\mathcal{E}_{7}(1) = \{(1,3),(1,12)\}$ and denote this edge operation by $o_{0}$. 
	In this way we obtain a derived colored graph $\mathcal{G}_{1}(\pi_{1})= \mathcal{G}(\pi,\mathcal{D}_{0},o_0)$ of $\mathcal{G}(\pi)$ with respect to $\mathcal{D}_{0}$ and $o_0$, that is depicted in Figure \ref{g:exa2}. 
	We proceed to compute a derived set $ \mathcal{D}_{c}(\mathcal{D}_{0}) = \{1,2,4,5,6,7,11\}$ of $\mathcal{D}_{0}$ in $\mathcal{G}_{1}(\pi_{1})$ as shown in Figure \ref{g:exa3} and denote this derived set by $\mathcal{D}_{1}$.
	Since $\mathcal{D}_{1} \neq V$ and $\mathcal{D}_{1} \neq \mathcal{D}_{0}$, the procedure will continue. 
	For the nodes $2 , 11 \in \mathcal{D}_{1}  $ in the graph $\mathcal{G}_{1}(\pi_{1})$, we have $ N_{V \setminus \mathcal{D}_{1} }(\{11\}) \subseteq N_{V \setminus \mathcal{D}_{1}}(2)$, and  for any $k \in N_{V \setminus \mathcal{D}_{1} }(\{11\})$, the two edges $(2,k)$ and $(11,k)$ have the same color. Thus we eliminate all the edges in $\mathcal{E}_{11}(2) = \{(2,12),(2,9)\}$ and denote this operation by $o_1$.  
	We then obtain a derived colored graph $\mathcal{G}_{2}(\pi_{2}) = \mathcal{G}_{1}(\pi_{1},\mathcal{D}_{1},o_1)$ of $\mathcal{G}_{1}(\pi_{1})$ with respect to $\mathcal{D}_{1}$ and $o_1$, and $\mathcal{G}_{2}(\pi_{2})$ is depicted in Figure  \ref{g:exa4}.  
	We then compute a derived set $\mathcal{D}_{c}(\mathcal{D}_{1})$ of $\mathcal{D}_{1}$ in $\mathcal{G}_{2}(\pi_{2})$ as shown in Figure \ref{g:exa5}. This derived set is denoted by  $\mathcal{D}_{2}$ and turns out to be equal to the original vertex set $V$.
		Thus we obtain that an edge-operations-color-change derived set $\mathcal{D}_{ec}(V_L)$ is equal to $V$, and conclude that $(\mathcal{G}(\pi);V_L)$ is controllable.	
\end{ex}

\section{Conclusion}

In this paper we have studied strong structural controllability of leader/follower networks. 
In contrast to existing work, in which the nonzero off-diagonal entries of matrices in the qualitative class are completely independent, in this paper we have studied the general case that there are equality constraints among these entries, in the sense that a priori given entries in the system matrix are restricted to take arbitrary but identical nonzero values.
This has been formalized using the concept of colored graph and by introducing the new concept of colored strong structural controllability.
In order to obtain conditions under which colored strong structural controllability holds for a given leader-follower system, we have introduced a new color change rule and a new concept of zero forcing set. 
These have been used to formulate a sufficient condition for controllability of the colored graph with a given leader set. We have shown that this condition is not necessary, by giving an example of a colored strong structurally controllable colored graph and leader set for which our sufficient condition is not satisified.

Motivated by this example, we have proceeded to establish the concept of elementary edge operations on colored graphs. It has been shown that these edge operations preserve colored strong structural controllability. Based on these elementary edge operations and the color change rule, a second sufficient graph theoretic condition for colored strong structural controllability has been provided. 

Finally, we have established a condition for colored strong structural controllability in terms of the new notion of edge-operations-color-change derived set. This derived set is obtained from the original leader set by applying edge operations and the color change rule sequentially in alternating manner. This iterative procedure has been illustrated by means of a concrete example.

The main new ideas of this paper are a new color change rule, and the concept of elementary edge operations for colored directed graphs. We have established several conditions for colored strong structural controllability using these new concepts. The conditions that we provided are not necessary, and finding necessary and sufficient conditions is still an open problem. Another open problem is to establish methods to characterize strong structural controllability for the case that given entries in the system matrices satisfy linear relations (instead of requiring them to take identical values). For {\em weak} structural controllability this was studied in \cite{LM2017}. 

In this paper we have focused on finding graph-theoretic conditions rather than providing suitable algorithms, see e.g. \cite{WRS2014}. Establishing an efficient algorithm to check colored strong structural controllability could also be a future research problem.   
Finally, other system-theoretic concepts like strong targeted controllability \cite{MCT2015, WCT2017} and identifiability \cite{WTC2018} for systems defined on colored graphs are possible research directions for the future.

\appendix[]
\noindent{\sc I: Proof of Lemma~\ref{l:BSC}.}
\medskip

\begin{IEEEproof}
By the Hautus test \cite{H1969},  $(\mathcal{G}(W); V_{L})$ is controllable if and only if $[A - \lambda I ~ B]$ has full row rank for all $A \in \mathcal{Q}_{W}(\mathcal{G})$ and all $\lambda \in \mathbb{C}$ with $B = B(V;V_{L})$ given by \eqref{e:inputmatrix}.
	Let $V = \{1,2, \ldots,n\}.$
	
	We first prove the `if' part. Suppose that $V_{L}$ is a balancing set for $\mathcal{G}(W)$.	
	Without loss of generality, we may assume that there is a chronological list of zero extensions
	\[(C_{1} \xrightarrow{z} Y_{1}, C_{2} \xrightarrow{z} Y_{2}, \ldots, C_{s} \xrightarrow{z} Y_{s}),\]
	where, for $r = 1, 2, \ldots, s$, $C_{r}$ represents the current set of zero vertices  before the $r$th zero extension and $Y_{r} \subseteq V \setminus C_{r}$, and  $C_{s} \cup Y_{s} = V$.
    Assign variables $x_{1},x_2,\ldots,x_{n}$ to every vertex in $V$, with $x_{i} = 0$ if $i \in C_{r}$ and $x_i$  undetermined otherwise. 
    To every vertex $j \in C_{r}$, we then assign a balance equation given by \eqref{e:BE}.
    By definition of the zero extension rule, we have the following implications
	\begin{equation}\label{e:BSBE1}
	x^{T}_{V \setminus C_{i}}W_{V \setminus C_{i},C_{i}} = 0 \Rightarrow x_{Y_{i}}^{T} = 0 \mbox{ for } i = 1,2, \ldots, s.
	\end{equation}  
	For any $A \in \mathcal{Q}_{W}(\mathcal{G})$ and $\lambda \in \mathbb{C}$, there exists a diagonal matrix  $D \in \mathbb{C}^{n \times n}$ such that $A - \lambda I = W + D$. 
	It then follows immediately that  \[(A-\lambda I)_{V \setminus C_{i},C_{i}} = W_{V \setminus C_{i},C_{i}} \mbox{ for } i = 1,2, \ldots, s.\]
	Recalling \eqref{e:BSBE1}, we have that 
    \[	x^{T}_{V \setminus C_{i}}(A-\lambda I)_{V \setminus C_{i},C_{i}} = 0 \Rightarrow x_{Y_{i}}^{T} = 0 \mbox{ for } i = 1,2, \ldots, s.\]
 Since $x^{T}B= 0 \Rightarrow x^{T}_{V_{L}} = 0$ and $ V_{L} \cup (\bigcup_{j=1}^{s}{Y}_{j}) = V $, we then have that
	\[x^{T}[A-\lambda I~ B]= 0 \Rightarrow x^{T} = 0,\]
	which implies that $[A-\lambda I~ B]$ has full row rank. 
	Since the $A$ and $\lambda$ are arbitrary, $(\mathcal{G}(W); V_{L})$ is controllable. 
	Thus we have proved the 'if' part.
	
	To prove the converse, suppose that $(\mathcal{G}(W); V_{L})$ is controllable while $V_{L}$ is not a balancing set. 
	It follows immediately that  $[A - \lambda I ~ B]$ has full row rank for all $A \in \mathcal{Q}_{W}(\mathcal{G})$ and all $\lambda \in \mathbb{C}$, with $B = B(V;V_{L})$ given by \eqref{e:inputmatrix}, and the derived set $D = \mathcal{D}_{z}(V_{L})$ is not equal to $V$.  
Again assign variables $x_{i}$ to the vertices $i \in V$ such that $x_{i} = 0$ if $i \in D$ and $x_i$ is undetermined otherwise.  Let $D' = V \setminus D$.	
	 By definition of the zero extension rule, we conclude that there exists a vector $x$ such that $x_{D} = 0$, $x_{D'} \neq 0$ and $x^T W = 0$, where $x_{D}$ and $x_{D'}$ are the sub-vectors  corresponding to the components in $D$ and $D'$, respectively. Recalling that $V_{L} \subseteq D$ , it follows  that
	 $x^T [W ~B] = 0$.
	  This implies that the matrix $[W~B]$ does not have full row rank.
	 Thus we have reached a contradiction and the proof is completed.
\end{IEEEproof}
\bigskip

\noindent {\sc II: example of non-uniqueness of derived sets.}
\medskip

\begin{ex} \label{ex:counterexample}
Consider the colored graph $\mathcal{G}(\pi) = (V,E,\pi)$  depicted in Figure \ref{g:ce1}. Take as coloring set $C = \{1,2,3,4,5\}$.  Consider the colored bipartite graph $G = (X,Y,E_{XY},\pi_{XY})$ associated with $X = \{1,2,3,4\}$ and $Y = \{6,7,8,9\}$ as is depicted in Figure \ref{g:ce2}. 
    It can be shown that there exists exactly one equivalence class of perfect matchings in $G$ with nonzero signature. Since $X \subset C$ and $Y = N_{V \setminus C}(X)$, we have that $X \xrightarrow{c} Y$. After applying this force we arrive at the derived set $\mathcal{D}_1(C) =  V$.
    
On the other hand, obviously $X_1 \xrightarrow{c} Y_1$, with $X_1 = \{5\}$ and $Y_1 = \{6\}$. After applying this force, no other forces are possible. Indeed, it can be verified that there does not exist a subset of $\{1, 2, 3, 4, 5, 6\}$ that forces any subset of $\{7, 8, 9\}$. In this way we arrive at the derived set $\mathcal{D}_2(C) = \{1,2, 3,4,5, 6\}$.

We conclude that there exist two different derived sets in $\mathcal{G}(\pi)$ with coloring set $C$. Thus we have found an example for the non-uniqueness of derived sets for a given colored graph and coloring set. 

\begin{figure}[h!]
	\centering
	\begin{subfigure}{0.4\textwidth}
		\centering
		\begin{tikzpicture}[scale=0.3]
		\tikzset{VertexStyle1/.style = {shape = circle,
				ball color = white!100!black,
				text = black,
				inner sep = 2pt,
				outer sep = 0pt,
				minimum size = 10 pt},
			edge/.style={->,> = latex', text = black}
		}
		\tikzset{VertexStyle2/.style = {shape = circle,
				ball color = black!80!yellow,
				text = white,
				inner sep = 2pt,
				outer sep = 0pt,
				minimum size = 10 pt}}
		\node[VertexStyle2](1) at (-2,2) {$1$};
		\node[VertexStyle2](2) at (10,3) {$2$};
		\node[VertexStyle2](3) at (5,-1) {$3$};
		\node[VertexStyle2](4) at (-2,-4) {$4$};
		\node[VertexStyle2](5) at (-1,7) {$5$};
		\node[VertexStyle1](6) at (-10,5) {$6$};
		\node[VertexStyle1](7) at (-9,-3) {$7$};
		\node[VertexStyle1](8) at (3,-7) {$8$};
		\node[VertexStyle1](9) at (1,0) {$9$};
		\Edge[ style = {->,> = latex'},color=red, label = $c_{1}$,labelstyle={inner sep=0pt}](1)(6);
		\Edge[ style = {->,> = latex'},color=red, label = $c_{1}$,labelstyle={inner sep=0pt}](1)(7);
		\Edge[ style = {->,> = latex'},color=red, label = $c_{1}$,labelstyle={inner sep=0pt}](1)(9);
		\Edge[ style = {->,> = latex'},color=green, label = $c_{2}$,labelstyle={inner sep=0pt}](2)(6);
		\Edge[ style = {->,> = latex'},color=green, label = $c_{2}$,labelstyle={inner sep=0pt}](2)(9);
		\Edge[ style = {->,> = latex'},color=blue, label = $c_{3}$,labelstyle={inner sep=0pt}](2)(8);
		\Edge[ style = {->,> = latex'},color=orange, label = $c_{5}$,labelstyle={inner sep=0pt}](3)(9);
		\Edge[ style = {->,> = latex'},color=blue, label = $c_{3}$,labelstyle={inner sep=0pt}](3)(8);
		\Edge[ style = {->,> = latex',pos = 0.3},color=green, label = $c_{2}$,labelstyle={inner sep=0pt}](4)(7);
		\Edge[ style = {->,> = latex'},color=black, label = $c_{4}$,labelstyle={inner sep=0pt}](4)(8);
		\Edge[ style = {->,> = latex'},color=orange, label = $c_{5}$,labelstyle={inner sep=0pt}](4)(9);
		\Edge[ style = {->,> = latex'},color=red, label = $c_{1}$,labelstyle={inner sep=0pt}](5)(6);
		\Edge[ style = {->,> = latex'},color=red, label = $c_{1}$,labelstyle={inner sep=0pt}](7)(6);
		\Edge[ style = {->,> = latex'},color=red, label = $c_{1}$,labelstyle={inner sep=0pt}](9)(8);
		\Edge[ style = {->,> = latex'},color=green, label = $c_{2}$,labelstyle={inner sep=0pt}](7)(9);
		\end{tikzpicture}
		\caption{A colored directed graph $\mathcal{G}(\pi)$.}
		\label{g:ce1}
	\end{subfigure}
	\medskip
	
	\begin{subfigure}{0.4\textwidth}
	\centering
	\begin{tikzpicture}[scale=0.9]
	\tikzset{VertexStyle1/.style = {shape = circle,
			ball color = white!100!black,
			text = black,
			inner sep = 2pt,
			outer sep = 0pt,
			minimum size = 10 pt},
		edge/.style={->,> = latex', text = black}
	}
	\tikzset{VertexStyle2/.style = {shape = circle,
			ball color = black!80!yellow,
			text = white,
			inner sep = 2pt,
			outer sep = 0pt,
			minimum size = 10 pt}}
	\node[VertexStyle2](1) at (-2.5,2.5) {$1$};
	\node[VertexStyle2](2) at (-2.5,0.8) {$2$};
	\node[VertexStyle2](3) at (-2.5,-0.8) {$3$};
	\node[VertexStyle2](4) at (-2.5,-2.5) {$4$};
	\node[VertexStyle1](6) at (2.5,2.5) {$6$};
	\node[VertexStyle1](7) at (2.5,0.8) {$7$};
	\node[VertexStyle1](8) at (2.5,-0.8) {$8$};
	\node[VertexStyle1](9) at (2.5,-2.5) {$9$};
	\Edge[ style = {> = latex',pos=0.3},color=red, label = $c_{1}$,labelstyle={inner sep=0pt}](1)(6);
	\Edge[ style = {> = latex',pos=0.3},color=red, label = $c_{1}$,labelstyle={inner sep=0pt}](1)(7);
	\Edge[ style = {> = latex',pos=0.3},color=red, label = $c_{1}$,labelstyle={inner sep=0pt}](1)(9);
	\Edge[ style = {> = latex',pos=0.3},color=green, label = $c_{2}$,labelstyle={inner sep=0pt}](2)(6);
	\Edge[ style = {->,> = latex',pos=0.3},color=green, label = $c_{2}$,labelstyle={inner sep=0pt}](2)(9);
	\Edge[ style = {> = latex',pos=0.3},color=blue, label = $c_{3}$,labelstyle={inner sep=0pt}](2)(8);
	\Edge[ style = {> = latex',pos=0.7},color=orange, label = $c_{5}$,labelstyle={inner sep=0pt}](3)(9);
	\Edge[ style = {> = latex',pos=0.3},color=blue, label = $c_{3}$,labelstyle={inner sep=0pt}](3)(8);
	\Edge[ style = {> = latex',pos = 0.15},color=green, label = $c_{2}$,labelstyle={inner sep=0pt}](4)(7);
	\Edge[ style = {> = latex',pos=0.3},color=black, label = $c_{4}$,labelstyle={inner sep=0pt}](4)(8);
	\Edge[ style = {> = latex',pos=0.3},color=orange, label = $c_{5}$,labelstyle={inner sep=0pt}](4)(9);
	\end{tikzpicture}
	\caption{Colored bipartite graph  $G$.}
	\label{g:ce2}
\end{subfigure}
\caption{Example of non-uniqueness of derived sets.}
\end{figure}
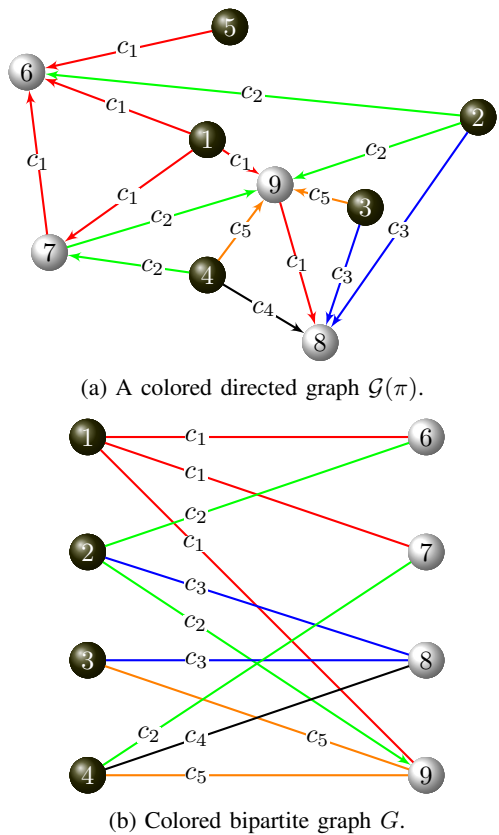
\end{ex}

%
%
%
%
%
%
%
%

\ifCLASSOPTIONcaptionsoff
  \newpage
\fi



%
%
%
%
\bibliographystyle{IEEEtran}
\bibliography{myieee2}

%

%
%
%




\end{document}